\documentclass[11pt]{amsart}
\usepackage[top=1.2in, bottom=1.5in, left=1.5in, right=1.5in]{geometry}

\usepackage{amsfonts, amssymb, amsmath, enumerate}

\numberwithin{equation}{section}

\DeclareMathOperator*{\im}{Im}

\DeclareMathOperator*{\Det}{det}

\def \N {\mathbb{N}}

\def \R {\mathbb{R}}

\def \e {\varepsilon}

\def \vol {{\rm vol}}
\def \Id {{\rm Id}}

\def \etc {,\ldots,}

\newcommand{\norm}[1]{\left \| #1 \right \|}
\newcommand{\pr}[2]{\left \langle {#1} , {#2} \right \rangle}

\newtheorem{theorem}{Theorem}[section]
\newtheorem{proposition}[theorem]{Proposition}
\newtheorem{corollary}[theorem]{Corollary}

\theoremstyle{remark}
\newtheorem{remark}[theorem]{Remark}

\begin{document}

\title[Non-central sections of the simplex, the cross-polytope and the cube]{Non-central sections of the simplex, the cross-polytope and the cube}

\author{Hermann K\"onig}
\address[Hermann K\"onig]{Mathematisches Seminar\\
Universit\"at Kiel\\
24098 Kiel, Germany
}
\email{hkoenig@math.uni-kiel.de}


\begin{abstract}
We determine the maximal hyperplane sections of the regular $n$-simplex, if the distance of the hyperplane to the centroid is fairly large, i.e. larger than the distance of the centroid to the midpoint of edges. Similar results for the $n$-cube and the $l_1^n$-ball were obtained by Moody, Stone, Zach and Zvavitch and by Liu and Tkocz. The maximal hyperplanes in these three cases are perpendicular to the vectors from the centroid to the vertices. For smaller distances -in a well-defined range- we show that these hyperplane sections are at least locally maximal. We also determine the hyperplane sections of the simplex, the cross-polytope and the cube which have maximal perimeter, i.e. maximal volume intersection with the boundary of the convex body.
\end{abstract}

\maketitle

\section{Introduction and main results}

The study of sections of convex bodies is a very active area with applications in functional analysis, probability and computer science. It is difficult to find the hyperplane sections of specific convex bodies which have maximal or minimal volume, even for classical bodies. In a well-known paper, Ball \cite{B} found the maximal hyperplane section of the $n$-cube. Meyer and Pajor \cite{MP} determined the maximal hyperplane sections of the cross-polytope which was extended by Koldobsky \cite{K}, Theorem 7.7, to the unit balls of $l_p^n$, $0 < p < 2$. Webb \cite{W} identified the maximal sections of the simplex through the centroid, using techniques of Ball. These results concern {\bf central} hyperplane sections through the origin in the case of symmetric bodies or through the centroid in the case of general convex bodies. \\

For non-central sections, not too many results are known. Moody, Stone, Zach and Zvavitch \cite{MSZZ} showed that the maximal hyperplane sections of the $n$-cube for large distances to the origin are the ones perpendicular to the main diagonals. Recently, Liu and Tkocz \cite{LT} proved a corresponding result for the cross-polytope, i.e. the $l_1^n$-unit ball, if the distance to the origin is larger than $\frac 1 {\sqrt2}$. In both cases the distance of the hyperplane to the origin has to be bigger than the distance from zero to the midpoint of the edges. \\

Motivated by these results, we first consider the corresponding problem for the regular $n$-simplex, also for large distances. For smaller distances, we prove  that the parallel hyperplanes are at least locally maximal. We show locally extremal results for the cross-polytope and the cube, too. Further, we prove that the same hyperplanes at large distances to the center have maximal perimeter, i.e. maximal intersection with the boundary of the convex body. The problem of the maximal perimeter of convex bodies originated with a question of Pelczy\'nski on the maximal perimeter of cubic sections. \\

To formulate our results precisely, we introduce some definitions and notations. Let $a \in S^{n-1} \subset \R^n$, i.e. $\norm {a} _2 =1$ and $t \in \R$. By
$$H_t(a) := \{ x \in \R^n \ | \ \pr{a}{x} = t \}  =  \{ t a\} + a^\perp $$
we denote the hyperplane perpendicular to $a$ at distance $t$ to the origin. Let $K \subset \R^n$ be a convex body with centroid in the origin. Then
$$A(a,t) = A_K(a,t) := \vol_{n-1}(H_t(a) \cap K) = \vol_{n-1}(\{ x \in K \ | \ \pr{x}{a} =t \}) $$
is the {\it parallel section function} of $K$ and
$$P(a,t) = P_K(a,t) := \vol_{n-2}(H_t(a) \cap \partial K) = \vol_{n-2}(\{ x \in \partial K \ | \ \pr{x}{a} =t \}) $$
is the {\it perimeter function} of $K$, where $\partial K$ denotes the topological boundary of $K$. For $t \ne 0$, we have non-central sections, and it is our aim to determine $\max_{a \in S^{n-1}} A(a,t)$ and $\max_{a \in S^{n-1}} P(a,t)$ for some polytopes and large $t$. If the body is clear from the context, we will omit the index $K$. As for convex bodies, we only study the unit $n$-cube $K = Q^n := [-\frac 1 2, \frac 1 2]^n$, the cross-polytope
$K=B(l_1^n) = \{ x \in \R^n \ | \ \norm {x} _1 \le 1 \}$ and the regular $n$-simplex $K= \Delta^n$ for which we use its representation in $\R^{n+1}$
$$\Delta^n := \{ x \in \R^{n+1} \ | \ x = (x_j)_{j=1}^{n+1} , \ x_j \ge 0 , \ \sum_{j=1}^{n+1} x_j =1 \} . $$
Since the centroid $c = \frac 1 {n+1} (1 \etc 1) \in \R^{n+1}$ of $\Delta^n$ is not in the origin when $\Delta^n$ is imbedded in $\R^{n+1}$ as above, we will assume in the definition of $A$ and $P$ that $a \in S^n \subset \R^{n+1}$ and $\sum_{j=1}^{n+1} a_j = 0$, so that $\pr{a}{c} = 0$. Then the hyperplane $H_0(a)$ passes through the centroid and the parallel hyperplane $H_t(a) = \{ x \in \R^{n+1} \ | \ \pr{a}{x} = t \}$ has distance $t$ to the centroid. We will assume in various results below that $t > d = d(n)$ is sufficiently large, where $d(n)$ is the distance of the centroid to the midpoint of edges of the convex body $K$. In this case, the hyperplane $H_t(a)$ will separate one vertex from all others, if it intersects $K$ at all. For $K = Q^n$, $d(n) = \frac {\sqrt{n-1}} 2$, and the range for $t$ is
$\frac {\sqrt{n-1}} 2 < t \le \frac {\sqrt n} 2$, cf. \cite{MSZZ}, for the cross-polytope it is $d(n) = \frac 1 {\sqrt 2} < t \le 1$, cf. \cite{LT} and for the simplex $d(n) = \sqrt{ \frac {n-1}{2(n+1)} } < t \le \sqrt { \frac n {n+1} }$. Note that $\sqrt { \frac n {n+1} }$ is the distance of the vertices
$e_j=(0 \etc 0, \underbrace{1}_j , 0 \etc 0)$ of $\Delta^n$ to the centroid $c$; it is the maximum coordinates of $a \in S^n$ with $\sum_{j=1}^{n+1} a_j = 0$ can attain, and in this case all other coordinates are equal to $-\frac 1 {\sqrt{n(n+1)}}$. The side-length of the simplex $\Delta^n$ is $\sqrt 2$, its height $\sqrt {\frac{n+1} n}$ and its volume $\frac{\sqrt{n+1}}{n!}$. The volume of the cross-polytope is $\vol_n(B(l_1^n)) = \frac {2^n}{n!}$. We can now state our main results. \\

\begin{theorem}\label{th1}
Let $n \ge 3$ and $K = \Delta^n \subset \R^{n+1}$ be the regular $n$-simplex and assume that $t \in \R$ satisfies $d(n) = \sqrt{ \frac {n-1}{2(n+1)} } < t \le \sqrt { \frac n {n+1} }$. \\
Let $a^{(n)} := (\sqrt{\frac n {n+1}},-\frac 1 {\sqrt{n(n+1)}} \etc -\frac 1 {\sqrt{n(n+1)}}) \in S^n$. Then $H_t(a^{(n)})$ is a maximal hyperplane section of $\Delta^n$ at distance $t$, i.e.
$$A(a,t) \le A(a^{(n)},t) = \frac {\sqrt{n+1}}{(n-1)!} \ (\frac n {n+1})^{n/2} \ (\sqrt{ \frac n {n+1}} -t )^{n-1}  $$
for all $a \in S^n$ with $\sum_{j=1}^{n+1} a_j = 0$. Up to permutations of coordinates of $a^{(n)}$, $H_t(a^{(n)})$ is the unique maximal hyperplane at distance $t$. For $n=2$, the same statement holds if $\frac 5 4 \frac 1 {\sqrt 6} \le t \le \sqrt {\frac 2 3}$. It is false for
$d(2) = \frac 1 {\sqrt 6} < t < \frac 5 4 \frac 1 {\sqrt 6}$.
\end{theorem}

\remark Since $a^{(n)}$ is a unit vector in the direction from the centroid $c$ to the vertex $e_1 = (1,0 \etc 0) \in \R^{n+1}$ of $\Delta^n$, the hyperplane $H_t(a)$ is parallel to the face $\{ x \in \Delta^n \ | \ x_1 = 0 \}$. The only other extremal hyperplanes are parallel to the other faces. \\

As a consequence, $V(a,t) := \vol_n(\{ x \in \Delta^n \ | \ \pr{a}{x} \ge t \})$ is also maximal for $a = a^{(n)}$, if $t > d(n)$. \\

If $t$ is smaller than $d(n)$, $a^{(n)}$ is still a local extremum of $A( . ,t)$:

\begin{proposition}\label{prop1}
Let $n \ge 3$ and $K = \Delta^n$ and $- \frac 1 {\sqrt{n(n+1)}} < t \le d(n) = \sqrt{ \frac {n-1}{2(n+1)} }$.
Let $c(n):= \frac{2n+1}{n(n+2)} \sqrt{\frac n {n+1}}$. Then \\
(a) If $c(n) < t \le d(n)$, $a^{(n)}$ is a local maximum of $A( . , t)$. \\
(b) If $- \frac 1 {\sqrt{n(n+1)}} < t < c(n)$, $a^{(n)}$ is a local minimum of $A( . ,t)$.
\end{proposition}

\remark For $n \ge 3$, $c(n) < d(n)$ and $\frac 2 {n+2} < c(n) < \frac 2 {n+1}$. Hence for $t> \frac 2 {n+1}$, $a^{(n)}$ is a local maximum and for $t < \frac 2 {n+2}$ a local minimum of $A( . ,t)$. Of course, the size of the neighborhood $\mathcal{U}$ of $a^{(n)}$ in which $a^{(n)}$ is a local extremum depends on $t$. The requirement is that $a=(a_1 \etc a_{n+1}) \in \mathcal{U}$ still satisfies $a_1 > t > a_2 \etc a_{n+1}$, just as the coordinates of $a^{(n)}$ do. In particular, for $t=0$, central sections of the simplex parallel to the faces are locally minimal, as has been shown already by Dirksen \cite{D}, Theorem 1.3. It has been claimed by Filliman \cite{F} that these sections are globally minimal, but no proof has been published. \\

For the perimeter of simplex sections we have

\begin{theorem}\label{th2}
Let $n \ge 5$, $K = \Delta^n$ and $d(n) = \sqrt{ \frac {n-1}{2(n+1)} } < t \le \sqrt { \frac n {n+1} }$. Then $H_t(a^{(n)})$ is a hyperplane section of the simplex at distance $t$ with maximal perimeter,
$$P(a,t) \le P(a^{(n)},t) = \frac{n \sqrt{n-1} }{(n-2)!} \ (\frac n {n+1})^{(n-2)/2} \ (\sqrt{\frac n {n+1}} -t)^{n-2} $$
for all $a \in S^n$ with $\sum_{j=1}^{n+1} a_j = 0$. Up to permutations of the coordinates of $a^{(n)}$, $H_t(a^{(n)})$ is the unique maximal hyperplane. For $n=4$, the same result holds, if $\sqrt{\frac 3 2} \sqrt{\frac 3 {10}} < t \le \sqrt{\frac 4 5}$.
\end{theorem}
The case $n=4$, $\sqrt{\frac 3 {10}} < t < \sqrt{\frac 3 2} \sqrt{\frac 3 {10}}$ is open.

\begin{proposition}\label{prop2}
For $n \ge 4$ and $t \in \R$ with $\frac{3n+2}{n(n+2)} \sqrt{\frac n{n+1} } < t \le d(n) = \sqrt{ \frac {n-1}{2(n+1)} }$, $a^{(n)}$ is a local maximum of the perimeter function $P( . , t)$ of the simplex.
\end{proposition}

\vspace{0,5cm}
As mentioned, Liu and Tkocz \cite{LT}, Theorem 4, determined the maximal hyperplane sections of the cross-polytope $K = B(l_1^n)$ at distance $t$, $\frac 1 {\sqrt 2} < t \le 1$. They showed that for these $t$, all $a \in S^{n-1} \subset \R^n$ with $n \ge 3$
\begin{equation}\label{eq1}
A(a,t) \le A(e_1,t) = \frac {2^{n-1}}{(n-1)!} (1-t)^{n-1} \ .
\end{equation}
For $n=2$, one needs $\frac 3 4 < t \le 1$ for this result to be true. The maximal hyperplanes are unique up to permutations of coordinates and reflections $x \to -x$. As for local extrema, we can show

\begin{proposition}\label{prop3}
Let $n \ge 3$, $K=B(l_1^n)$ be the cross-polytope and $0 < t \le \frac 1 {\sqrt 2}$, $t \ne \frac 3 {n+2}$. Then the axis-parallel hyperplane $H_t(e_1)$ provides a local extremum of the parallel section function $A( . ,t)$, namely \\
(a) a local maximum if $\frac 3 {n+2} < t \le \frac 1 {\sqrt2}$ and \\
(b) a local minimum if $0 < t < \frac 3 {n+2}$. \\
For $n=2$, we have a local minimum for $0 < t < \frac 3 4$.
\end{proposition}

\remark Again the size of the neighborhood, where $e_1$ yields a local extremum of $A( .,t)$, depends on $t$. For $n=2$, $e_1$ actually is the absolute minimum of $A(.,t)$, if $1 - \frac 1 {\sqrt 2} < t < \frac 1 {\sqrt 2}$, cf. \cite{LT}, Theorem 2. By Meyer and Pajor \cite{MP}, $e_1$ is the absolute maximum of the central section function $A(.,0)$ of the cross-polytope, although for $t>0$ close to zero, $e_1$ is a local minimum of $A(.,t)$ (in a very small neighborhood of $e_1$). These results show that $e_1$ is not a (relative or absolute) maximum of $A(.,t)$ for all $t$. \\

For the perimeter of $l_1^n$-sections, we have

\begin{theorem}\label{th3}
Let $n \ge 4$, $K=(l_1^n)$ be the cross-polytope and $\frac 1 {\sqrt 2} < t \le 1$. Then $H_t(e_1)$ is a maximal perimeter hyperplane section of the cross-polytope: For all $a \in S^{n-1}$,
$$P(a,t) \le P(e_1,t) = \frac{\sqrt{n-1}}{(n-2)!} \ 2^{n-1} \ (1-t)^{n-2} \ . $$
For $n=3$, the same is true at least for $\frac 4 5 < t \le 1$.
\end{theorem}

\remark The numerical evidence is that the result holds for $n=3$ and $t$ with $\frac 1 {\sqrt 2} < t \le \frac 4 5$, too. The proof of the previous result also yields

\begin{corollary}\label{cor1}
Let $n \ge 6$, $\frac 4 n < t \le \frac 1 {\sqrt 2}$. Then $e_1$ is a local maximum of the perimeter function $P( ., t)$ of the cross-polytope.
\end{corollary}

As mentioned, Moody, Stone, Zach and Zvavitch \cite{MSZZ} found the maximal hyperplane sections of the unit volume $n$-cube for large distances $t$, a result which started this line of investigations. Their result is

\begin{theorem}\label{th4}
Let $n \ge 3$, $t \in \R$ with $\frac {\sqrt{n-1}} 2 < t \le \frac {\sqrt n} 2$ and $a^{[n]} := \frac 1 {\sqrt n} (1 \etc 1) \in S^{n-1}$. Then the hyperplane $H_t(a^{[n]})$ perpendicular to the main diagonal $a^{[n]}$ is maximal at distance $t$: For all $a \in S^{n-1}$,
$$A(a,t) \le A(a^{[n]},t) = \frac {n^{n/2}}{(n-1)!} \ (\frac {\sqrt n} 2 -t )^{n-1} . $$
For $n=2$, one needs $\frac 3 8 \sqrt 2 < t \le \frac 1 {\sqrt 2}$; for $\frac 1 2 < t < \frac 3 8 \sqrt 2$ the statement is not true.
\end{theorem}

We provide an alternative proof of this result for $n \ge 3$ which also yields information on local maxima if
$\frac{n-2}{2 \sqrt n} < t \le \frac {\sqrt {n-1}} 2$. The narrow range of values for $t$ results from the fact that intersecting hyperplanes should separate one vertex of $Q^n$, say $e=\frac 1 2 (1 \etc 1)$, from all others, in particular from $e^{[i]} = \frac 1 2 ( 1, \etc ,1 , \underbrace{-1}_i, 1, \etc 1)$ and that
$\pr{a^{[n]}}{e} = \frac {\sqrt n} 2$ while $\pr{a^{[n]}}{e^{[i]}} = \frac {n-2} {2 \sqrt n} $. We have

\begin{proposition}\label{prop4}
Let $n \ge 5$, $K=Q^n$ and $\frac{n-2}{2 \sqrt n} < t \le \frac {\sqrt {n-1}} 2$. Then $a^{[n]}$ is a local maximum of $A( . , t)$.
\end{proposition}

More precisely, for any such $t$, let \\
$\mathcal{U} := \{ a \in S^{n-1} \ | \  \sum_{j=1}^n a_j > 2 t > \sum_{j=1}^n a_j - 2 a_i  \text{  for all   } i= 1 \etc n \}$.
Then $\mathcal{U}$ is a neighborhood of $a^{[n]}$ and $A(a^{[n]},t) = \max_{a \in \mathcal{U}} A(a,t)$.
Note that $a^{[n]} \in \mathcal{U}$ since $\sqrt n - \frac 2 {\sqrt n} = \frac {n-2}{\sqrt n} < 2t$. \\
The result is not true for $n=2$; $a^{[2]}$ yields the absolute {\it minimum} of $A( . , t)$ if $\frac 1 {\sqrt 2} - \frac 1 2 < t < \frac 1 2$ and a relative minimum if $0 < t < \frac 1 {\sqrt 2} - \frac 1 2$, cf. K\"onig, Koldobsky \cite{KK1}, Proposition 5. For $n=3$, $a^{[3]}$ gives the absolute maximum if $0.642 < t < \frac {\sqrt 3} 2$, but the absolute minimum if $0.32 < t < \frac 1 2$, and a relative minimum if $\frac 1 {2 \sqrt 3} < t < 0.32$, cf. Proposition 7 of \cite{KK1}. For $n=4$, $a^{[4]}$ is a relative maximum if $\frac 5 8 < t \le 1$ and a relative minimum if $ \frac 1 2 < t < \frac 5 8$, as the proof of Proposition \ref{prop4} will show. Therefore Proposition \ref{prop4} not valid in the full range of $t$-values when $n=2, 3, 4$. \\

As for the maximal perimeter of the $n$-cube, we have

\begin{theorem}\label{th5}
Let $n \ge 4$, $K=Q^n$ and $\frac{\sqrt {n-1}} 2 < t \le \frac {\sqrt n} 2$. Then the hyperplane section $H_t(a^{[n]})$ of the cube orthogonal to the main diagonal  has maximal perimeter among all hyperplanes of distance $t$ to the origin, i.e.
$$P(a,t) \le P(a^{[n]},t) = \frac{\sqrt{n-1}}{(n-2)!} \ n^{n/2} \ (\frac 1 2 \sqrt n - t)^{n-2}  $$
for all $a \in S^{n-1}$.
The result is true for $n=3$, too, if $0.725 < t \le \frac {\sqrt 3} 2$.
\end{theorem}
Numerical evidence indicates that the result is also true for $n=3$ if $\frac 1 {\sqrt 2} < t < 0.725$. As for local maxima we have

\begin{proposition}\label{prop5}
Let $n \ge 6$, $\frac{n-2}{\sqrt n} < t \le \frac {\sqrt{n-1}} 2$. Then the main diagonal $a^{[n]}$ is a local maximum of the perimeter function $P( . ,t)$ of the cube $Q^n$.
\end{proposition}
The result holds for a more restricted range of $t$-values also in dimensions $n= 3, 4, 5$: For $n=5$, $a^{[n]}$ is a local maximum of $P(.,t)$ if
$\frac 7 {4 \sqrt 5} < t \le 1$, for $n=4$ if $\frac 3 4 < t < \frac {\sqrt 3} 2$ and for $n=3$, if $0.722 \simeq \frac 5 {4 \sqrt 3} < t \le \frac 1 {\sqrt 2}$. However, for $n=3$ and $\frac 1 {\sqrt 3} < t < \frac {11}{30} \sqrt 3 \simeq 0.635$, $a^{[3]}$ is actually a local minimum of $P(.,t)$, as we will see. \\

\remark Let us mention that for central cubic sections, $t=0$, it was shown by K\"onig and Koldobsky \cite{KK} that the vector
$a(2):= \frac 1 {\sqrt 2}(1,1,0 \etc 0) \in \R^n$, $n \ge 3$ yields the maximal perimeter, i.e. we have for all $a \in S^{n-1}$
that $P(a,0) \le P(a(2),0)$. This corresponds to Ball's result \cite{B} that $A(a,0) \le A(a(2),0)$.

\remark Not much is known on lower estimates of $A( . ,t)$ and $P( . ,t)$ for non-central sections. K\"onig and Rudelson \cite{KR} proved for the cube
$K=Q^n$ and $t=\frac 1 2$ that $A(a,t) \ge \frac 1 {17}$ for all $n \in \N$ and $a \in S^{n-1}$. Note that for $t > \frac 1 2$, $A(a,t)$ may be zero, e.g. for $a = e_1$. It is open whether the minimum of $A( . , \frac 1 2)$ is attained for the main diagonal $a^{[n]}$ or not. For $n=2$ and $n=3$ this is true, cf. K\"onig and Koldobsky \cite {KK1}.

\section{Hyperplane sections of the simplex}

Let $\Delta^n := \{ x \in \R^{n+1} \ | \ x = (x_j)_{j=1}^{n+1} , \ x_j \ge 0 , \ \sum_{j=1}^{n+1} x_j =1 \} $ be the $n$-simplex. It has centroid
$c = \frac 1 {n+1} (1 \etc 1) \in \R^{n+1}$ and vertices $e_j$, $j \in \{1  \etc n+1\}$. For $a \in \R^{n+1}$ with $\norm {a} _2 =1$, $\sum_{j=1}^{n+1} a_j = 0$, let $|a_1| = \max_{1 \le j \le n+1} |a_j|$. We may assume that $a_1 >0$. Let $t > d(n) = \sqrt{ \frac{n-1}{2(n+1)} }$. If
$H_t(a) = \{ x \in \R^{n+1} \ | \ \pr{a}{x} = t \}$ intersects $\Delta^n$ non-trivially, $a_1 = \pr{a}{e_1} \ge t$. All other coordinates of $a$ are less than $t$, since if e.g. $a_2 \ge t$,
$\frac{a_1+a_2} 2 \ge t > d(n) = \norm{\frac{e_1+e_2} 2 -c}_2 \ge \pr{a}{\frac{e_1+e_2} 2 -c} = \frac{a_1+a_2} 2 $
would yield a contradiction. Hence $H_t(a)$ separates the vertex $e_1$ from all others. Since for $a_1=t$, $A(a,t) = \vol_{n-1}(H_t(a) \cap \Delta^n) = 0$, we may assume in the following that $a_1 > t > a_2 \etc a_{n+1}$. \\

\begin{proposition}\label{prop2.1}
Let $a \in \R^{n+1}$, $\norm{a}_2 =1$ with $\sum_{j=1}^{n+1} a_j =0$. Let $t \in \R$ be such that $a_1 > t > a_2 \etc a_{n+1}$. Then
\begin{equation}\label{eq2.1}
A(a,t) = \frac {\sqrt{n+1}}{(n-1)!} \ \prod_{j=2}^{n+1} \frac 1 {a_1-a_j} \ (a_1 -t)^{n-1} \ .
\end{equation}
\end{proposition}

\remark Since $|a_1| \le \sqrt{ \frac n {n+1} }$, also $|t| \le \sqrt{ \frac n {n+1} }$. Formula \eqref{eq2.1} follows from a modification of Corollary 2.4 of Dirksen \cite{D}, but we will give a simple direct geometric proof. \\

\begin{proof}
The hyperplane $H_t(a)$ intersects the edges from $e_1$ to $e_j$ in points $v_j = s_j e_1 + (1-s_j) e_j$, $j \in\{2 \etc n+1\}$, with $\pr{a}{v_j} = t$. Let $P$ be the $n$-dimensional pyramid spanned by the vectors $v_j-e_1 = (1-s_j) (e_j-e_1)$, i.e. let $P$ be the convex hull of $e_1$ and $H_t(a) \cap \Delta^n$. The pyramid has height $a_1-t$ since $\pr{a}{e_1-v_j} = a_1 - t$. Therefore
$$\vol_n(P) = \frac 1 n \vol_{n-1}(H_t(a) \cap \Delta^n) (a_1-t) = \frac {a_1-t} n \ A(a,t)  \ . $$
Let $B$ be the $n \times (n+1)$-matrix $B = ((1-s_j) (e_j-e_1))_{j=2 \etc n+1}$. Then $\vol_n(P) = \frac 1 {n!} \sqrt{\det(B^* B)}$,
$$\Det(B^* B) = \prod_{j=2}^{n+1} (1-s_j)^2 \Det (\Id + J) \ , $$
where $J$ is the matrix of all ones. Since $\Id + J$ has one eigenvalue $n+1$, with all other eigenvalues being $1$, $\Det(\Id +J) = n+1$. Therefore
$$\vol_n(P) = \frac{\sqrt{n+1}}{n!} \ \prod_{j=2}^{n+1} (1-s_j) = \frac{a_1-t} n \ A(a,t) \ . $$
Since $t = \pr{a}{v_j} = s_j a_1 + (1-s_j) a_j$, $a_1-t = (1-s_j)(a_1-a_j)$, $1-s_j = \frac{a_1-t}{a_1-a_j}$. Therefore
$$A(a,t) =  \frac{\sqrt{n+1}}{(n-1)!} \prod_{j=2}^{n+1} \frac 1 {a_1-a_j} \ (a_1-t)^{n-1} \ . $$
\end{proof}

{\it Proof of Theorem \ref{th1}.} \\
(1) Let $\mathcal{S} := \{ x \in \R^{n+1} \ | \ \norm{x} _2 =1, \ \sum_{j=1}^{n+1} x_j =0 \}$ be the surface of constraints. To find $\max A(a,t)$ or
$\max \ln A(a,t)$ for $a \in \mathcal{S}$, we first determine the critical points of $\ln A( . , t)$ by Lagrangian multipliers. We also assume the restrictions $a_1 > t > a_2 \etc a_{n+1}$ so that Proposition \ref{prop2.1} applies. Looking at
$$F(a,\lambda,\mu) := (n-1) \ln(a_1 - t) - \sum_{j=2}^{n+1} \ln(a_1 - a_j) + \frac \lambda 2 (\sum_{j=1}^{n+1} a_j^2 - 1) + \mu (\sum_{j=1}^{n+1} a_j ) \ , $$
we have to solve the equations
$$\frac {\partial F}{\partial a_1} = \frac {n-1}{a_1-t} - \sum_{j=2}^{n+1} \frac 1 {a_1-a_j} + \lambda a_1 + \mu = 0 $$
$$\frac {\partial F}{\partial a_j} =  \frac 1 {a_1-a_j} + \lambda a_j + \mu = 0 \ , \ j=2 \etc n+1 \ $$
in $\mathcal{S}$. Adding these $n+1$ equations, we get $\mu = - \frac{n-1}{n+1} \frac 1 {a_1-t}$. Inserting this into the first equation yields
$$\frac {n(n-1)}{n+1} \frac 1 {a_1-t} - \sum_{j=2}^{n+1} \frac 1 {a_1-a_j} + \lambda a_1 = 0 \ . $$
Multiplying this by $a_1$ and the $n$ equations for $j = 2 \etc n+1$ by $a_j$ and adding these gives
$$n = \sum_{j=2}^{n+1} \frac{a_1-a_j}{a_1-a_j} = \frac {n(n-1)}{n+1} \frac {a_1} {a_1-t} + \lambda + \mu \sum_{j=2}^{n+1} a_j \ , \ \sum_{j=2}^{n+1} a_j = -a_1 \ , $$
$$\lambda = n - (n-1) \frac {a_1}{a_1-t} = - \frac {nt -a_1}{a_1-t} \ . $$
From $\frac {\partial F}{\partial a_j} = 0$ we get $-(\lambda a_j + \mu) = \frac 1 {a_1-a_j}$ for all $j=2 \etc n+1$. Inserting the values of $\lambda$ and $\mu$, the coordinates $a_j$ of a critical point $a \in \mathcal{S}$ have to be solutions of the quadratic equation
$$a_j^2 - p \ a_j -q = 0 \quad , \quad p:= a_1 - \frac{n-1}{n+1} \frac 1 {n t -a_1} \quad , \quad q:= \frac{(n+1) \ t - 2 a_1}{(n+1)(n t -a_1)} \ .$$
We now assume that $t \ge \frac 2 {\sqrt{n(n+1)}}$ and claim that $t^2-p \ t - q \le 0$. Calculation shows
$$t^2 - p \ t -q = -\frac {a_1-t}{(n+1)(nt-a_1)} \ \Big( \ n(n+1) t^2 -(n+1) a_1 -2 \ \Big) \ . $$
The function $\phi(t) := n(n+1) t^2 -(n+1) a_1 -2$ is non-negative if $t \ge t_+$, where
$t_+ = \frac{a_1}{2n} + \sqrt{ \frac{a_1^2}{4 n^2} + \frac 2 {n(n+1)} }$. Since $a_1 \le \sqrt{ \frac n {n+1} }$, $t_+ \le \frac 2 {\sqrt{n(n+1)}}$. Thus by assumption on $t$, $t \ge t_+$ and hence $t^2 - p \ t -q \le 0$. This implies that the larger solution $a_{j,+} = \frac p 2 + \sqrt {\frac {p^2} 4 + q}$ of the quadratic equation for $a_j$ satisfies $a_{j,+} \ge t$ which contradicts our assumption that $a_j < t$ for $j =2 \etc n+1$. Therefore only one solution, namely
$a_{j,-} = \frac p 2 - \sqrt {\frac {p^2} 4 + q}$ for $a_j$ is allowed, and all $a_j$ for $j \ge 2$ have to be equal,
$a_2 = \dots a_{n+1} = \sqrt{ \frac {1-a_1^2} n }$. Since $\sum_{j=1}^{n+1} a_j = 0$ is required, too, $a_1 = \sqrt {\frac n {n+1} }$,
$a_j = - \frac 1 {\sqrt{n(n+1)}}$, so that $a = a^{(n)} := (\sqrt{\frac n {n+1}},-\frac 1 {\sqrt{n(n+1)}} \etc -\frac 1 {\sqrt{n(n+1)}})$ is the only possible critical value under our restrictions. Clearly, $a^{(n)}$ is a solves the Lagrange equations for appropriate $\lambda$ and $\mu$, and thus is a critical point. By assumption in Theorem \ref{th1}, $ t > \sqrt{ \frac{n-1}{2(n+1)} } > \frac 2 {\sqrt{n(n+1)} } $ is satisfied for any $n \ge 4$. For $n=3$,
$\frac 1 2 \le t_+ \le \frac 1 {\sqrt 3}$. Hence for $n=3$ and $t > \frac 1 {\sqrt 3}$, the same conclusion $a = a^{(3)}$ holds. The cases $n=3$,
$\frac 1 2 \le t \le \frac 1 {\sqrt 3}$ and $n=2$ will be considered separately. \\

(2) We now show that $a^{(n)}$ is a relative maximum of $A( . ,t)$. Since it was the only critical point in the allowed range for $a \in \mathcal{S}$, it will then be  the absolute maximum. Let $b : [-\e_0,\e_0] \to \mathcal{S}$ be a smooth, locally not constant $C^2$-function with $b(0) = a^{(n)}$, $b_1 > t > b_j$, $j=2, \dots , n+1$. We show that for $f(a,t) := \prod_{j=2}^{n+1} \frac 1 {a_1-a_j} \ (a_1-t)^{n-1}$, $\ln f \circ b$ has a maximum in $0$ for all such curves, fixing $t$. Then $\ln f$ and $f$ have a local maximum in $b(0)=a^{(n)}$. Clearly $(\ln f \circ b)'(0) = 0$ because $a^{(n)}$ is a critical value. Since $b$ maps into $\mathcal{S}$, we have $\sum_{j=1}^{n+1} b_j = 0$, $\sum_{j=1}^{n+1} b_j^2 = 1$ which implies
$$\sum_{j=1}^{n+1} b_j' = 0 \ , \ \sum_{j=1}^{n+1} b_j b_j' = 0 \ , \ \sum_{j=1}^{n+1} b_j'' = 0 \ , \ \sum_{j=1}^{n+1} ((b_j')^2+b_j b_j'') = 0 \ . $$
In particular, $b_1'(0) = - \sum_{j=2}^{n+1} b_j'(0)$, $\sqrt{\frac n {n+1}} \ b_1'(0) = - \frac 1 {\sqrt{n(n+1)}} \sum_{j=2}^{n+1} b_j'(0)$, hence $b_1'(0)=0$ and $\sum_{j=2}^{n+1} b_j'(0)=0$. Moreover, $b_1''(0) = - \sum_{j=2}^{n+1} b_j''(0)$ and
$\sqrt{\frac n {n+1}} \ b_1''(0) = - \frac 1 {\sqrt{n(n+1)}} \sum_{j=2}^{n+1} b_j''(0) + \sum_{j=2}^{n+1} b_j'(0)^2$, yielding
$b_1''(0) = -\sqrt { \frac n {n+1} } \sum_{j=2}^{n+1} b_j'(0)^2 < 0$. \\

Calculation shows
$$(\ln f \circ b)'' = \frac{n-1}{b_1-t} \ b_1'' - \frac{n-1}{(b_1-t)^2} \ (b_1')^2 + \sum_{j=2}^{n+1} \frac{(b_1'-b_j')^2}{(b_1-b_j)^2} - \sum_{j=2}^{n+1} \frac{b_1''-b_j''}{b_1-b_j} \ , $$
$$(\ln f \circ b)''(0) = \frac{n-1}{\sqrt{\frac n {n+1}}-t} \ b_1''(0) + \frac n {n+1} \sum_{j=2}^{n+1} b_j'(0)^2 - \sqrt{\frac n {n+1}} \ \sum_{j=2}^{n+1} (b_1''(0)-b_j''(0)) \ . $$
Since  $\sum_{j=2}^{n+1} (b_1''(0) - b_j''(0)) = n \ b_1''(0) - \sum_{j=2}^{n+1} b_j''(0) = (n+1) \ b_1''(0)$ \\
and $\sum_{j=2}^{n+1} b_j'(0)^2 = - \sqrt{\frac{n+1}n } \ b_1''(0)$, we find
\begin{equation}\label{eq2.2}
(\ln f \circ b)''(0) = \Big(\frac{n-1}{\sqrt{\frac n {n+1}}-t} - (n+2) \sqrt{ \frac n {n+1} } \Big) \ b_1''(0) \ .
\end{equation}
Now $\frac{n-1}{\sqrt{\frac n {n+1}}-t} - (n+2) \sqrt{ \frac n {n+1} } > 0$ if and only if $t > \sqrt{\frac n {n+1}} \frac{2n+1}{n(n+2)} =:c(n)$.
By assumption, $t > d(n) = \sqrt{\frac{n-1}{2(n+1)}}$. Since $d(n) \ge c(n)$ for all $n \ge 3$ and $b_1''(0) <0$, $(\ln f \circ b)''(0) <0$. Therefore $a^{(n)}$ is a maximum of $f( . , t)$ in $\mathcal{S}$. If $n=2$, we have a relative maximum if $t > \frac 5 4 \frac 1 {\sqrt 6} =c(2)$. This is assumed in Theorem \ref{th1} if $n=2$. \\

(3) For $n=3$ and $\frac 1 2 < t < \frac 1 {\sqrt 3}$, there are possibly two values for $a_j$, $j=2,3,4$. Suppose that $a_2=a_3$. Then $a_1^2+2a_2^2+a_4^2=1$,
$a_1+2a_2+a_4=0$,
$$a_2 = \frac 1 6 (\pm \sqrt{6-8 a_1^2} -2 a_1) \quad , \quad a_4 = \frac 1 6 (\mp 2 \sqrt{6-8 a_1^2} -2 a_1) \ . $$
For $f(a,t) := \frac{(a_1-t)^2}{(a_1-a_2)^2 (a_1-a_4)}$, calculation shows that, when $a_2$ and $a_4$ are expressed as functions of $a_1$, the derivative of
$\ln f (a,t)$ with respect to $a_1$ has the same sign as $\phi(a_1,t) := 30 a_1 t -3 -10 a_1^2 +(3 t - a_1) \sqrt{6 - 8 a_1^2}$. We claim that $\phi$ is positive. Since $\phi$ is increasing in $t$, the minimum occurs for $t = \frac 1 2$,
$\phi(a_1,\frac 1 2) = 15 a_1 -3 -10 a_1^2 +(\frac 3 2 - a_1) \sqrt{6-8 a_1^2} \ge 2 >0$. Hence the maximum of $f$ occurs for the maximal possible value of $a_1$, namely $\sqrt{\frac 3 4}$. \\

(4) For $n=2$, $A(a,t) = \sqrt 3 \ \frac{a_1-t}{(a_1-a_2)(a_1-a_3)}$ with $a_1^2+a_2^2+a_3^2=1$, $a_3=-a_1-a_2$, implying $a_1^2+a_2^2+a_1 a_2 - \frac 1 2 = 0$. Thus $(a_1-a_2)(a_1-a_3) = (a_1-a_2)(2 a_1+a_2) = 2 a_1^2 - a_1 a_2 -a_2^2 = 3 a_1^2 - \frac 1 2$, $A(a,t) = \sqrt 3 \ \frac{a_1-t}{3 a_1^2 - \frac 1 2}$. Since
$\frac d {d a_1} A(a,t) = \sqrt 3 \ \frac{6 a_1 t -3 a_1^2 -\frac 1 2}{(3 a_1^2- \frac 1 2)^2} \ge 0$ if and only if $t \ge \frac 1 2 (a_1+\frac 1 {6 a_1})$,
$A( .,t)$ attains its maximum at $a^{(2)} = \frac 1 {\sqrt 6} (2,-1,-1)$ if $t \ge \max_{a_1 \in [\frac 1 {\sqrt 6}, \frac 2 {\sqrt 6}]} \frac 1 2 (a_1 + \frac 1 {6 a_1}) = \frac 5 4 \frac 1 {\sqrt 6} \simeq 0.51$, as assumed in Theorem \ref{th1}. \\
For $\frac 1 {\sqrt 6} < t < \frac 5 4 \frac 1 {\sqrt  6}$, the maximum of $A( . , t)$ is attained for $a \ne a^{(2)}$, namely for
$$a_1 = t + \sqrt{t^2-\frac 1 6} < \frac 2 {\sqrt 6} \quad , \quad a_{2,3} = -\frac t 2 - \frac 1 2 \sqrt{t^2 - \frac 1 6} \pm \sqrt{\frac 5 8 - \frac 3 2 t^2- \frac 3 2 t \sqrt{t^2-\frac 1 6} } \ . $$    \hfill $\Box$

{\it Proof of Proposition \ref{prop1}.} \\
Let $\mathcal{S}$ be as in the previous proof and $- \frac 1 {\sqrt{n(n+1)}} < t \le d(n) = \sqrt{\frac{n-1}{2(n+1)}}$. Define $\mathcal{U} := \{ x \in \mathcal{S} \ | \ x_1 > t > x_2 \etc x_{n+1} \}$. Then $a^{(n)} \in \mathcal{U}$ and $\mathcal{U}$ is an open neighborhood of $a^{(n)}$ in $\mathcal{S}$. The formula given in Proposition \ref{prop2.1} is valid for all $a \in \mathcal{U}$. Part (1) of the previous proof shows that
$a^{(n)}$ is a critical point of $A( . , t)$ in $\mathcal{U} \subset \mathcal{S}$, though maybe not the only one if $t < \frac 2 {\sqrt{n(n+1)}}$. Part (2) with equation \eqref{eq2.2} shows that $a^{(n)}$ is a relative maximum if $t > c(n) = \sqrt{\frac n {n+1}} \frac{2n+1}{n(n+2)}$ since then
$(\ln f \circ b)''(0) < 0$. For $t < c(n)$, $(\ln f \circ b)''(0) > 0$, and $a^{(n)}$ is a relative minimum. \\

Note that the neighborhood $\mathcal{U}$ of $a^{(n)}$ is rather small if $t$ is close to $- \frac 1 {\sqrt{n(n+1)}}$. We cannot claim that $a^{(n)}$ yields the absolute maximum of $A( . ,t)$ in $\mathcal{S}$, even if there is only one critical point, since for these small values of $t$ the formula for $A(a,t)$ will be different from the one in Proposition \ref{prop2.1}, if $H_t(a)$ no longer separates one vertex from all others.    \hfill $\Box$

\remark
We sketch a different proof of Theorem \ref{th1}. One has to show for any $a \in \mathcal{S}$, which attains the maximum of $A( . , t)$, that $a_2 = \dots = a_{n+1}$, under the assumption that $a_1 > t > a_2 \etc a_{n+1}$. If e.g. $a_2 \ne a_3$, replace $a_1, a_2, a_3$ by $a_1' = a_1 + \e$, $a_2' = a_3' = \frac {a_2+a_3- \e} 2$ and let $a_j' = a_j$ for $j>3$. Then $\sum_{j=1}^{n+1} a_j' = 0$,  and $\sum_{j=1}^{n+1} a_j'^2 =1$ holds if
$\frac 1 2 (a_2-a_3)^2 = \e (2 a_1 -a_2-a_3) + \frac 3 2 \e^2$. This equation admits a positive solution $\e > 0$ since $a_2 \ne a_3$. We claim that
$A(a,t) < A(a',t)$: This will imply that all coordinates except for the first have to be equal in the maximal situation. Now $A(a,t) < A(a',t)$ is equivalent to
$$(1 + \frac {\e}{a_1-t})^{n-1} > \frac{(a_1-\frac{a_2+a_3} 2 + \frac{3 \e} 2)^2}{(a_1-a_2)(a_1-a_3)} \ \prod_{j=4}^{n+1} (1+ \frac{\e}{a_1-a_j}) \ . $$
As in the proof of Proposition \ref{prop2.1}, let $v_j = s_j e_1 + (1-s_j) e_j \in H_t(a) \cap span (e_1,e_j)$,
$s_j a_1 + (1-s_j) a_j = \pr{a}{v_j} > t > \norm{\frac{e_1+e_j} 2 - c}_2 \ge \pr{a}{\frac{e_1+e_j} 2 - c} = \frac{a_1+a_j} 2 $, hence
$(s_j- \frac 1 2)(a_1-a_j) > 0$, $s_j > \frac 1 2$. Since $1-s_j = \frac{a_1-t}{a_1-a_j}$, $a_1-t \le \frac 1 2 (a_1-a_j)$ and $\frac{\e}{a_1-t} \ge \frac 1 2 \frac{\e}{a_1-a_j}$. Hence it suffices to show
$$(1+ \frac{\e}{a_1-t}) (1+ \frac{\e}{2(a_1-t) + \e})^{n-2} > \frac{(a_1-\frac{a_2+a_3} 2 + \frac{3 \e} 2)^2}{(a_1-a_2)(a_1-a_3)} =
1+ \frac{(a_2-a_3)^2}{(a_1-a_2)(a_1-a_3)} \ .$$
Since $f(x):=(1+x)(1+\frac x {2+x})^{n-2} > 1 + \frac n 2 x$ for all $x >0$, choosing $x = \frac{\e}{a_1-t}$, it suffices to prove that
$$ \frac n 2 \frac{\e}{a_1-t} \ge \frac{(a_2-a_3)^2}{(a_1-a_2)(a_1-a_3)} \ . $$
It follows from the equation for $\e$ that $\e \ge \frac 1 3 \frac{(a_2-a_3)^2}{2 a_1 - a_2 - a_3}$. Using this and $a_2, \ a_3 < a_1$ as well as
$\frac{\e}{a_1-t} > 2 \frac{\e}{a_1-a_2}$, it follows that
$$ 3 \frac{\e}{a_1-t} \ge \frac{(a_2-a_3)^2}{(a_1-a_2)(a_1-a_3)} \ . $$
This proves the claim if $n \ge 6$ and also for $n=5$ if the estimate for $f$ is done a bit more carefully. The proof may be modified also in the situation of Proposition \ref{prop1} when $t < d(n)$ and $a_1 > t > a_2 \etc a_{n+1}$. Then there is $\alpha < 1$ such that $a_1-t < \alpha (a_1-a_j)$ for all $j=2 \etc n+1$ and the modified proof will rely on the estimate $\frac{\e}{a_1-t} \ge \frac 1 {\alpha} \frac{\e}{a_1-a_j}$.    \hfill $\Box$

\vspace{0,5cm}
To prove Theorem \ref{th2} on the maximal perimeter of simplex sections, we first derive a formula for the perimeter.

\begin{proposition}\label{prop2.2}
Let $a \in \R^{n+1}$, $\norm{a}_2 = 1$, $\sum_{j=1}^{n+1} a_j = 0$, $t \in \R$ be such that $a_1 > t > a_2 \etc a_{n+1}$. Then
\begin{equation}\label{eq2.3}
P(a,t) = \frac 1 {(n-2)!} \ \sum_{j=2}^{n+1} \sqrt{n - (n+1) \ a_j^2} \ \prod_{k=2,k \ne j}^{n+1} \frac 1 {a_1-a_k} \ (a_1-t)^{n-2} \ .
\end{equation}
\end{proposition}

\begin{proof}
By definition, $P(a,t) = \vol_{n-2} (\{ x \in \Delta^n \ | \ \pr{a}{x} = t \} \cap \partial \Delta^n)$. The boundary of $\Delta^n$ consists of the \ $n+1$ \ simplices $(\partial \Delta^n)_j = \{x \in \Delta^n \ | \ x_j = 0 \}$ of dimension $(n-1)$. The hyperplane $H_t(a)$ does not intersect $(\partial \Delta^n)_1$ since we have for $x_1=0$ that $\pr{a}{x} = \sum_{j=2}^{n+1} a_j x_j < t$ because $\sum_{j=2}^{n+1} x_j = 1$, $a_j < t$. \\
Consider $H_t(a) \cap (\partial \Delta^n)_{n+1}$, i.e. $x_{n+1}=0$. Denote $\tilde{x}=(x_1 \etc x_n) \in \Delta^{n-1}, \tilde{a}=(a_1 \etc a_n) \in \R^n$. Let $e:=(1 \etc 1) \in \R^n$ and put $\tilde{b} = \frac{\tilde{a} + \frac{a_{n+1}} n e}{\sqrt{1-\frac{n+1}n a_{n+1}^2}}$. Then $\sum_{j=1}^n \tilde{b}_j = 0$ since $\sum_{j=1}^n a_j = - a_{n+1}$ and
$\sum_{j=1}^n \tilde{b}_j^2 = 1$ since $\sum_{j=1}^n a_j^2 = 1- a_{n+1}^2$. Now \\
$$\pr{\tilde{a}}{\tilde{x}} = t \text{ \ is equivalent to \ } \pr{\tilde{b}}{\tilde{x}} = \frac{t+\frac{a_{n+1}} n}{\sqrt{1-\frac{n+1} n a_{n+1}^2}} =: t_{n+1} \ . $$
This defines a hyperplane section of the $(n-1)$-simplex $(\partial \Delta^n)_{n+1}$ with $\tilde{b}_1 > t_{n+1} > \tilde{b}_2 \etc \tilde{b}_n$. Thus by Proposition \ref{prop2.1} the perimeter part of $P(a,t)$ on $(\partial \Delta^n)_{n+1}$ is
\begin{align*}
P_{(n+1)}(a,t) &= \frac{\sqrt n}{(n-2)!} \ \prod_{j=2}^n \frac 1 {\tilde{b}_1-\tilde{b}_j} \ (\tilde{b}_1 - t_{n+1})^{n-2} \\
& = \frac{\sqrt n}{(n-2)!} \ \sqrt{1 - \frac{n+1} n \ a_{n+1}^2} \ \prod_{j=2}^n \frac 1 {a_1-a_j} \ (a_1 - t)^{n-2} \ .
\end{align*}
Summing over all boundary simplices $(\partial \Delta^n)_j$ for $j=2 \etc n+1$ yields the formula
$$P(a,t) = \frac 1 {(n-2)!} \ \sum_{j=2}^{n+1} \sqrt{n - (n+1) \ a_j^2} \ \prod_{k=2,k \ne j}^{n+1} \frac 1 {a_1-a_k} \ (a_1-t)^{n-2} \ . $$
\end{proof}

There is also a general formula for $P(a,t)$ when the hyperplane $H_t(a)$ does not separate one vertex from all others, namely
\begin{equation}\label{eq2.g}
P(a,t) = \frac 1 {(n-2)!} \sum_{j=1}^{n+1} \sqrt{ n - (n+1) \ a_j^2 } \ \frac 1 {2 \pi} \int_\R \prod_{k=1, k \ne j}^{n+1} \frac 1 {1+ i (a_k - t) s } \ ds \ . \end{equation}
This can be shown similarly as above, when the $(n-2)$-dimensional area of \\ $H_h(a) \cap (\partial \Delta^n)_j$ is calculated according to Dirksen \cite{D}, Theorem 1.1, Corollary 2.4 and the remarks there. \\

\vspace{0,5cm}
{\it Proof of Theorem \ref{th2}.} \\
(i) Since $t > d(n)$, $H_t(a)$ separates one vertex of $\Delta^n$, say $e_1$, from all others. Then $a_1 = \pr{a}{e_1} > t > \pr{a}{e_j} = a_j$ for all $j= 2 \etc n+1$. Hence by \eqref{eq2.3}
\begin{equation}\label{eq2.4}
P(a,t) = \frac 1 {(n-2)!} \ \sum_{j=2}^{n+1} \sqrt{n - (n+1) a_j^2} \ (a_1-a_j) \ \prod_{k=2}^{n+1} \frac 1 {a_1-a_k} \ (a_1-t)^{n-2} \ .
\end{equation}
Note that the product extends over all $k = 2 \etc n+1$ since we moved the factor $a_1-a_j$ to the weights.
Let $f(a):= \prod_{k=2}^{n+1} \frac 1 {a_1-a_j} \ (a_1-t)^{n-2}$. We claim that $f$, or equivalently $\ln f$, attains it maximum on
$\mathcal{S} := \{ x \in \R^{n+1} \ | \ \norm{x} _2 =1, \sum_{j=1}^{n+1} x_j =0 \}$ in $a^{(n)}$, and uniquely so. A modification of the proof of Theorem \ref{th1} works, except that the power $(n-1)$ is replaced by $(n-2)$ and thus the degree of homogeneity in the numerator and in the denominator is different. We solve the Lagrange equations
$\frac {\partial F}{\partial a_j} = 0$ for
$$F(a,\lambda,\mu) := (n-2) \ln(a_1 - t) - \sum_{j=2}^{n+1} \ln(a_1 - a_j) + \frac \lambda 2 (\sum_{j=1}^{n+1} a_j^2 - 1) + \mu (\sum_{j=1}^{n+1} a_j ) \ . $$
We find that
$$\lambda = - \frac{n-2}{n+1} \frac 1 {a_1-t} \quad , \quad \mu = - \frac{n t - 2 a_1}{a_1-t} $$
and the quadratic equation for any $a_j$, $ j \in \{2 \etc n+1 \}$ now reads
$$a_j^2 - p \ a_j - q = 0 \quad , \quad p = a_1 - \frac{n-2}{n+1} \frac 1 {n t - 2 a_1} \quad , \quad q = \frac{(n+1) t - 3 a_1}{(n+1)(n t -2 a_1)} \ . $$
Assume that $t > \frac 3 {\sqrt{n(n+1)}}$. Then $ t^2 - p \ t - q \le 0$ : We have
$$t^2 - p \ t - q = - \frac{a_1-t}{(n+1)(n t - 2 a_1)} \Big( n(n+1) t^2 - 2(n+1) a_1 t - 3 \Big) $$
and $\phi(t) := n(n+1) t^2 - 2 (n+1) a_1 t - 3$ is non-negative, if $t \ge t_+$, where
$t_+ = \frac{a_1} n + \sqrt{\frac{a_1^2}{n^2} + \frac 3 {n(n+1)}}$. Since $a_1 \le \sqrt{\frac n {n+1}}$, $t_+ \le \frac 3 {\sqrt{n(n+1)}}$,
$t^2 - p \ t -q \le 0$. Therefore $a_{j,+} \ge t$ which contradicts the assumption that $a_j < t$. Hence under the condition $a_j < t < a_1$, we have only one possible solution $a_j=a_{j,-}$, i.e. $a_2 = \dots = a_{n+1} = \frac{1-a_1^2} n$. Since also $\sum_{j=1}^{n+1} a_j = 0$, this again forces $a=(a_1 \etc a_{n+1})$ to be equal to $a^{(n)}$; $a^{(n)}$ is a critical point and the only critical point under the assumptions made. Note that $\frac 3 {\sqrt{n(n+1)}} < d(n) = \sqrt{\frac{n-1}{2(n+1)}}$ if $n \ge 5$. For $n=4$, there may be further critical points besides $a^{(4)}$, if $d(4) = \sqrt{\frac 3 {10}} < t < \frac 3 {\sqrt{20}}$ since then possibly $t_+ < t$. However, for $\frac 3 {\sqrt{20}} < t < \sqrt{\frac 4 5}$, $a^{(4)}$ is the only critical point of $f$ with $a_1>t>a_j$. \\

(ii) We show that $a^{(n)}$ is a relative maximum. Since it is the only critical point, it will be the absolute maximum of $f$ on $\mathcal{S}$. We proceed as in the proof of Theorem \ref{th1}, part (2). For any smooth non-constant function $b : [-\e_0,\e_0] \to \mathcal{S}$ with $b(0) = a^{(n)}$, $b_1 > t > b_2 \etc b_{n+1}$,  we find that, fixing $t$,
\begin{equation}\label{eq2.5}
(\ln f \circ b)''(0) = \Big ( \frac{n-2}{\sqrt{\frac n {n+1}} -t} - (n+2) \sqrt{\frac n {n+1}} \Big ) \ b_1''(0) \quad , \quad b_1''(0) < 0 .
\end{equation}
We need $\frac{n-2}{\sqrt{\frac n {n+1}} -t} > (n+2) \sqrt{\frac n {n+1}}$ to assure that $a^{(n)}$ is a relative maximum, meaning that
$t > \frac{3n+2}{n(n+2)} \sqrt{\frac n {n+1}}$ is needed. But by the assumption in Theorem \ref{th2}
$ t > d(n) = \sqrt{\frac{n-1}{2(n+1)}} > \frac{3n+2}{n(n+2)} \sqrt{\frac n {n+1}}$ holds for any $n \ge 5$. For $n=4$, we assumed $t > \frac 3 {\sqrt{20}}$, so $ t > \frac 7 {10} \sqrt{\frac 4 5}$ is also satisfied. \\

(iii) We now estimate $P(a,t)$ using \eqref{eq2.4}
\begin{equation}\label{eq2.6}
P(a,t) \le \frac 1 {(n-2)!} \ \Big (\sum_{j=2}^{n+1} \sqrt{n - (n+1) \ a_j^2} \ (a_1-a_j) \Big ) \ f(a^{(n)},t) \ ,
\end{equation}
$f(a^{(n)},t) = (\frac n {n+1})^{n/2} (\sqrt{\frac n {n+1}} -t )^{n-2}$. By Cauchy-Schwarz, using $a_1 \le \sqrt{\frac n {n+1}}$,
\begin{align*}
( \ \sum_{j=2}^{n+1} \sqrt{n - (n+1)  \ a_j^2} \ & (a_1-a_j) \ ) \le ( \ \sum_{j=2}^{n+1} (n-(n+1)a_j^2) \ )^{1/2} ( \ \sum_{j=2}^{n+1} (a_1-a_j)^2 \ )^{1/2} \\
&= ( \ n^2 - (n+1) (1-a_1^2) \ )^{1/2} ( \ \sum_{j=2}^{n+1} (a_1^2+a_j^2 - 2 a_1 a_j) \ )^{1/2} \\
&= ( \ n^2-n-1+(n+1)a_1^2 \ )^{1/2} ( \ na_1^2+(1-a_1^2)-2a_1(-a_1) \ )^{1/2} \\
&= ( \ n^2-n-1+(n+1)a_1^2 \ )^{1/2} ( \ 1+(n+1) a_1^2 \ )^{1/2} \\
&\le (n^2-1)^{1/2} (n+1)^{1/2} = \sqrt{n-1} \ (n+1) \ .
\end{align*}
Thus for all $a \in \mathcal{S}$
$$P(a,t) \le \frac{\sqrt{n-1}}{(n-2)!} \frac{n^{n/2}}{(n+1)^{n/2-1}} \ (\sqrt{\frac n {n+1}}-t )^{n-2} = P(a^{(n)},t) \ , $$
so that $a^{(n)}$ attains the maximum of $P( . ,t)$ in $\mathcal{S}$.        \hfill $\Box$

\vspace{0,5cm}
{\it Proof of Proposition \ref{prop2}.}  Let $n \ge 4$. By the previous proof, $a^{(n)}$ is a critical point of $f$ and actually a relative maximum since
$t > c(n) = \frac{3n+2}{n(n+2)} \sqrt{\frac n{n+1} }$ is assumed. The sum of weights in front of $f$ in \eqref{eq2.6} assume their absolute maximum in
$a^{(n)}$ as seen in part (iii) of the previous proof. Hence $P( . ,t)$ has a relative maximum in $a^{(n)}$. Note that $c(n) < \frac 3 {n+ 2/3}$. Hence for $t > \frac 3 {n+ 2/3}$, $a^{(n)}$ is a relative maximum of the perimeter function of the simplex.  \hfill $\Box$

\vspace{0,3cm}
\remark We did not investigate the maximal perimeter of the simplex in dimension 3 in detail, which might be done using the formula of Proposition \ref{prop2.2}. However, for $n=3$ and $t=d(3) = \frac 1 2$ or $t$ close to it, $t=\frac 1 2 + \delta$, $\delta$ small, $P(a,t) \le P(a^{(3)},t)$ is not true for all $a \in \mathcal{S}$:
Consider
$\bar{a} = (\frac 1 2 + 2 \e, \frac 1 2 - 2 \e, -\frac 1 2, -\frac 1 2)/\sqrt{1+8 \e^2} \in S^3$ for small $\e > 0$, $\sum_{j=1}^4 \bar{a}_j = 0$,
$\bar{a}_1 = \frac{\frac 1 2 + 2 \e}{\sqrt{1+8 \e^2}} > \frac 1 2 = d(3) > \bar{a}_2, \bar{a}_3, \bar{a}_4$. Let $Q = (\frac 1 2 + 2 \e, \frac 1 2 - 2 \e,0,0)$ be on the line from $e_1$ to $e_2$, $\bar{a}_1 > t := \norm{Q-c}_2 = \frac 1 2 \sqrt{1+8 \e^2} > \frac 1 2$. The perimeter formula \eqref{eq2.3} yields
$$P(a^{(3)},t) = (\frac 9 4 - \frac 3 4 \sqrt 6) - 3 \sqrt 6 \ \e^3 + O(\e^3) \ , $$
$$P(\bar{a},t) = \sqrt 2 - 2 \sqrt 2 \ \e + 6 \sqrt 2 \ \e^2 + O(\e^3) \ , $$
so that with $\frac 9 4 - \frac 3 4 \sqrt 6 \simeq 1.345 < 1.414 \simeq \sqrt 2$ we have $P(a^{(3)},t) < P(\bar{a},t)$ for $t = \frac 1 2 \sqrt{1+8 \e^2}$. Geometrically, $H_t(\bar{a}) \cap (\partial \Delta^n)$ is a distorted triangle, with two sides almost equal to $\frac 1 {\sqrt 2}$ and the third one of order $O(\e)$. \hfill $\Box$

\vspace{0,5cm}
As mentioned, Webb \cite{W} found the maximal {\it central} sections of the simplex. He showed that for all $a \in \mathcal{S}$,
$$A(a,0) \le A(a^{[2]},0) \quad , \quad a^{[2]} :=\frac 1 {\sqrt 2} (1,-1,0 \etc 0) \in S^n$$
The corresponding question whether the maximal central perimeter of the simplex is also attained by $a^{[2]}$ is open, i.e. whether $P(a,0) \le P(a^{[2]},0)$ holds for all $a \in \mathcal{S}$ or not. We are only able to prove an asymptotic result.

\begin{proposition}
Let $n \ge 3$, $a \in S^n \subset \R^{n+1}$ with $\sum_{j=1}^{n+1} a_j = 0$. Then
$$P(a,0) \le \frac{\sqrt{n-1}}{(n-2)!} \ (n+1) \ \sqrt{1 + \frac 1 n} \le (1 + \frac 1 n) \ P(a^{[2]},0) \ . $$
\end{proposition}

\begin{proof}
Let $n \ge 3$, $a \in S^n$, $\sum_{j=1}^{n+1} a_j = 0$. Then $a$ has at least two non-zero coordinates. By \eqref{eq2.g}
$$P(a,0) = \frac 1 {(n-2)!} \ \sum_{j=1}^{n+1} \sqrt{n - (n+1) a_j^2} \ D_j(a) \ , $$
$$D_j(a) := \frac 1 {2 \pi} \int_\R \prod_{k=1, k \ne j}^{n+1} \frac 1 {1+ i a_k s} \ ds \ . $$
Let $f(s):= \prod_{k=1}^{n+1} (1+i a_k s)$, $g_j(s):= 1 + i a_j s$. By perturbation, we may assume that all non-zero $a_k$'s occur only once in $f$. Then
$\frac 1 f$ has only simple poles at $s= \frac i {a_k}$ and
$$res_{\frac i {a_k}} ( \frac {g_j} f ) = \frac{ g_j(\frac i {a_k})}{f'(\frac i {a_k})} = \frac{1-\frac{a_j}{a_k}}{f'(\frac i {a_k})} = (1-\frac{a_j}{a_k}) \ res_{\frac i {a_k}} (\frac 1 f) \ . $$
Since $\frac {g_j} f$ is the integrand in $D_j(a)$, the residue theorem implies, using $\sum_{j=0}^{n+1} a_j = 0$,
\begin{align*}
\sum_{j=1}^{n+1} D_j(a) &= \sum_{j=1}^{n+1} \sum_{\im (\frac i {a_k}) > 0} (1- \frac{a_j}{a_k}) \ res_{\frac i {a_k}}(\frac 1 f) \\
&= \sum_{\im (\frac i {a_k}) > 0} \sum_{j=1}^{n+1} (1- \frac{a_j}{a_k}) \ res_{\frac i {a_k}}(\frac 1 f) \\
&= (n+1) \sum_{\im (\frac i {a_k}) > 0} res_{\frac i {a_k}}(\frac 1 f) = (n+1) \ \frac 1 {2 \pi} \int_\R \prod_{k=1}^{n+1} \frac 1 {1+i a_k s} \ ds \ .
\end{align*}
By Webb's result \cite{W} on maximal central sections of the simplex, we know that
$$\frac 1 {2 \pi} \int_\R \prod_{k=1}^{n+1} \frac 1 {1+i a_k s} \ ds \le \frac 1 {\sqrt 2} \ . $$
Therefore $\sum_{j=1}^{n+1} D_j(a) \le \frac{n+1}{\sqrt 2}$. Equation \eqref{eq2.g} then gives
$$P(a,0) \le \frac 1 {(n-2)!} \ \sqrt n \ \sum_{j=1}^{n+1} D_j(a) \le \frac{\sqrt n}{(n-2)!} \frac{n+1}{\sqrt 2} \ . $$

An application of the residue theorem to formula \eqref{eq2.g} yields for $a=a^{[2]}$ that
$$P(a^{[2]},0) = \frac{\sqrt{n-1}}{(n-2)!} \ ( \sqrt{n(n-1)} \ \frac 1 {\sqrt 2} + 1) \ . $$
Hence for all $a \in S^n$, $\sum_{j=1}^{n+1} a_j = 0$
$$P(a,0) \le (1 + \frac 1 n) \ P(a^{[2]},0) \ , $$
since $\frac{n+1}{n-1 + \sqrt{2 \frac{n-1} n}} \le 1 + \frac 1 n$ holds for all $n \ge 3$.
\end{proof}

\vspace{0,5cm}
It is not always true that the extrema for the parallel section function $A$ and the perimeter function $P$ are attained for the same vector.\\

{\bf Example.} Let $\Delta^3$ be the regular $3$-simplex. The {\it minimal} central section is attained for
$a^{(3)} = (\sqrt{\frac 3 4},-\frac 1 {2 \sqrt 3}, -\frac 1 {2 \sqrt 3}, -\frac 1 {2 \sqrt 3})$, $A(a,0) \ge A(a^{(3)},0)$ for all $a \in \mathcal{S}$, cf. Dirksen \cite{D}, Theorem 1.3. Now consider $\bar{a} = \frac 1 2 (1,-1,1,-1)$. Then $H_0(\bar{a}) \cap \Delta^3$ is a square of side-length $\frac 1 {\sqrt2}$ while  $H_0(a^{(3)}) \cap \Delta^3$ is a triangle of side-length $\frac 3 4 \sqrt 2$. Therefore
$$\frac 1 2 = A(\bar{a},0) > A(a^{(3)},0) = \frac 9 {32} \sqrt 3 \simeq 0.487 $$
while
$$2.828 \simeq 2 \sqrt 2 = P(\bar{a},0) < P(a^{(3)},0) = \frac 9 4 \sqrt 2 \simeq 3.182 \ . $$
Central hyperplane sections of the $3$-simplex are triangles or rectangles. The minimal central perimeter of triangle sections is attained for $a^{(3)}$, the minimal one for rectangles by $\bar{a}$. Thus $P(\bar{a},0) = 2 \sqrt 2$ is the minimal central perimeter of the $3$-simplex. \\
Similar examples with alternating $\pm 1$-sequences $\bar{a}$ can be given in dimensions $n=5 ,7$ with $P(\bar{a},0) < P(a^{(n)},0)$. However for odd dimensions $n \ge 9$ and similar alternating $\pm 1$-vectors $\bar{a}$ we have $P(\bar{a},0) > P(a^{(n)},0)$.

\section{Hyperplane sections of the cross-polytope}

We now study hyperplane sections of the cross-polytope, i.e. of the unit ball of $l_1^n$, $K = \{ x \in \R^n \ | \ \norm{x}_1 \le 1 \}$. $K$ has volume $\frac{2^n}{n!}$, side-length $\sqrt 2$ and vertices $\pm e_j$. The distance of the midpoints of edges to zero is $\frac 1 {\sqrt 2}$.
Let $a \in S^{n-1} \subset \R^n$ and $d = \frac 1 {\sqrt 2} < t \le 1$. If $H_t(a) \cap K$ is non-void, $H_t(a)$ separates one vertex, say $e_1$ from all others. Then $a_1 = \pr{a}{e_1} > t > \pr{a}{\pm e_j} = \pm a_j$. Hence we may and will assume in this section that $a_1 > t > a_j \ge 0$ holds for all $j = 2 \etc n$. \\
Liu and Tkocz proved in \cite{LT} the following formula for the parallel section function
\begin{equation}\label{eq3.1}
A(a,t) = \frac{2^{n-1}}{(n-1)!} \ a_1^{n-2} \ \prod_{j=2}^n \frac 1 {a_1^2 -a_j^2} \ (a_1 -t)^{n-1} \ .
\end{equation}
Their proof also works if $0 < t \le \frac 1 {\sqrt 2}$ as long as $a_1 > t > a_j$, $j =2 \etc n$ is satisfied which we assume. Similar as in Proposition \ref{prop2.1}, the formula follows from calculating the $n$-dimensional volume of the chopped-off part $P = \{ x \in K \ | \ \pr{a}{x} \ge t \}$ of the cross-polytope. Using \eqref{eq3.1} they proved that for $\frac 1 {\sqrt 2} < t \le 1$ and all $a \in S^{n-1}$ that 
\begin{equation}\label{eq3.2}
A(a,t) \le A(e_1,t) = \frac{2^{n-1}}{(n-1)!} \ (1-t)^{n-1} \ .
\end{equation}

\vspace{0,5cm}
{\it Proof of Proposition \ref{prop3}.} \\
(1) Let $a \in S^{n-1}$, $t>0$ with $a_1 > t > a_2 \etc a_n \ge 0$ and \\
$$f(a,t):= a_1^{n-2} \ \prod_{j=2}^n \frac 1 {a_1^2 -a_j^2} \ (a_1 -t)^{n-1} \ . $$
To find relative extrema of $f(.,t)$, we may as well find them for $\ln f(.,t)$. Thus, using Lagrange multipliers, we investigate
$$F(a,\lambda) = (n-2) \ln a_1 - \sum_{j=2}^n \ln (a_1^2-a_j^2) + (n-1) \ln (a_1-t) + \frac {\lambda} 2 (\sum_{j=1}^n a_j^2 -1) $$
for critical points on $S^{n-1}$. We have to solve
$$\frac {\partial F}{\partial a_1} = \frac {n-2}{a_1} - \sum_{j=2}^n \frac {2 a_1} {a_1^2-a_j^2} + \frac{n-1}{a_1-t} + \lambda a_1  = 0 $$
$$\frac {\partial F}{\partial a_j} =  \frac {2 a_j} {a_1^2-a_j^2} + \lambda a_j = 0 \ , \ j= 2 \etc n \ . $$
Of course, $a= e_1$ is a critical point, with $\lambda = n - \frac 1 {1-t}$. If there are other critical points, some $a_j$'s have to be non-zero. Then
$\lambda = - \frac 2 {a_1^2-a_j^2}$, so all non-zero $a_j$'s have to be equal. Assume $m$ such $a_j$'s are non-zero,
$a_2 = \dots = a_{m+1} \ne 0 = a_{m+2} = \dots = a_n$, with $m \le n-1$ since $a_1 \ge t > a_j$. Then $a_2^2 = \frac{1-a_1^2} m$,
$a_1^2 - a_j^2 = \frac {(m+1) a_1^2 -1 } m$ for $j = 2 \etc m+1$ and $\lambda = -\frac{2 m}{(m+1) a_1^2 -1}$. Inserting these values into the first equation yields
$$\frac{2 m - n} {a_1} - \frac{2 m (m+1) a_1}{(m+1) a_1^2-1} + \frac{n-1}{a_1-t} = 0 \ . $$
Solving this for $t$ gives
\begin{equation}\label{eq3.2a}
t = \frac{(m+2) a_1^2 + 2 m - 1}{(m+1) n a_1^2 + 2 m - n} a_1 =: \Phi(n,m,a_1) \ .
\end{equation}
We require that $0 < t \le a_1$. This is equivalent to $\frac 1 {\sqrt{m+1}} \le a_1$. Then the denominator $N:=(m+1) n a_1^2 + 2 m - n \ge 2m > 0$ of $\Phi(n,m,a_1)$ is positive.
The function $\Phi$ is decreasing in $m$ for all $n$ and $a_1$ such that $N$ is positive, since
$\frac{\partial \Phi}{\partial m} = -\frac 2 {N^2} \ (n-1) a_1 (1-a_1^2) < 0$. Thus $\Phi(n,n-1,a_1) \le \Phi(n,m,a_1) \le \Phi(n,1,a_1)$ (as long as $N$ is positive). We have $\Phi(n,m,\frac 1 {\sqrt 2}) \le \Phi(n,1,\frac 1 {\sqrt 2}) = \frac 1 {\sqrt 2}$ and
$\frac{\partial \Phi}{\partial {a_1}} = -\frac 1 {N^2} \Big( n-2+4na_1^2(2-a_1^2)-12a_1^2 \Big) < 0$ for all $a_1 > \frac 1 {\sqrt 2}$. Hence for all
$a_1 > \frac 1 {\sqrt 2}$, $\Phi(n,m,a_1) \le \Phi(n,1,a_1) < \frac 1 {\sqrt 2}$. Therefore equation \eqref{eq3.2a} has no solution with
$a_1 \ge t > \frac 1 {\sqrt 2}$. This means that for $t > \frac 1 {\sqrt 2}$, $e_1$ is the only critical point of $f( . ,t)$ with $a_1 \ge t > a_j$. \\

The smallest function $\Phi(n,n-1,a_1)$ is decreasing in $a_1 \in {(\frac 1 {\sqrt n},1)}$ if $n \ge 6$ since
$\frac{\partial \Phi(n,n-1,a_1)}{\partial {a_1}} = - \frac 1 {N^2} \Big( (n a_1^2-1)( \ (2-a_1^2) n^2 + 6 - 7 n \ ) \Big)< 0$, with
$\Phi(n,n-1,1) = \frac 3 {n+2} =: M(n)$. In fact, for all $m=1 \etc n-1$, $\Phi(n,m,1) = \frac 3 {n+2}$. For $n=3, 4, 5$, $\Phi(n,n-1, .)$ has a minimum $M(n)$ in some $\bar{a_1} \in (\frac 1 {\sqrt 2}, 1)$ slightly smaller than $\frac 3 {n+2}$, namely $\frac 1 {\sqrt 3}, \ \frac 5 {32} \sqrt{10}, \ \frac 7 {75} \sqrt{21}$, respectively. Then for $0 < t < M(n)$, equation \eqref{eq3.2a} has no solution, so in this case $e_1$  is the only critical value of $f( . ,t)$ with
$a_1 > t >0$, too. If $M(n) < t < \frac 1 {\sqrt 2}$, there are further critical points of $f( . , t)$ besides $e_1$, of the form
$\bar{a} = (a_1, \underbrace{a_2 \etc a_2}_m , 0 \etc 0)$ for a suitable $1 \le m \le n-1$, which are found by solving \eqref{eq3.2a} for $a_1$ and putting $a_2 = \sqrt{\frac{1-a_1^2} m}$. \\

(2) Let $b : [-\e_0,\e_0] \to S^{n-1}$ be a smooth $C^2$-function with $b(0) = e_1$ and $b_1 > t > b_2 \etc b_n$. Put
$$f(a,t):=a_1^{n-2} \ \prod_{j=2}^n \frac 1 {a_1^2 -a_j^2} \ (a_1 -t)^{n-1} \ . $$
To show that $f$ or $\ln f$ has a local maximum or minimum in $e_1$, it suffices to prove this for all functions $\ln f \circ b$, $b$ of the above type. Since $\sum_{j=1}^n b_j^2 =1$, we have
$$ \sum_{j=1}^n b_j b_j' = 0 \quad , \quad \sum_{j=1}^n ((b_j')^2 + b_j b_j'') = 0 \ . $$
With $b(0) = e_1$ we find that $b_1'(0)=0$ and $b_1''(0) = - \sum_{j=2}^n b_j'(0)^2 < 0$. Calculation gives for
$\ln f \circ b = (n-2) \ln b_1 + (n-1) \ln (b_1 - t) - \sum_{j=2}^n \ln (b_1^2 -b_j^2)$ that
\begin{align*}
(\ln f \circ b)'' &= - \frac{n-2}{b_1^2} \ (b_1')^2 - \frac{n-1}{(b_1-t)^2} \ (b_1')^2 + \frac{n-2}{b_1} \ b_1'' + \frac{n-1}{b_1-t} \ b_1'' \\
& + \sum_{j=2}^n \frac{(2 b_1 b_1' - 2 b_j b_j')^2}{(b_1^2 - b_j^2)^2} - \sum_{j=2}^n \frac{2 b_1'^2 - 2 b_j'^2 + 2 b_1 b_1'' - 2 b_j b_j''}{b_1^2-b_j^2} \ .
\end{align*}
In zero, most terms vanish, yielding
$$(\ln f \circ b)''(0) = \Big( (n-2) + \frac{n-1}{1-t} \Big) \ b_1''(0) + 2 \sum_{j=2}^n b_j'(0)^2 - \sum_{j=2}^n (2 \ b_1''(0)) \ . $$
Since $2 \sum_{j=2}^n b_j'(0)^2 = -2 \ b_1''(0)$, this implies
$$(\ln f \circ b)''(0) = \Big( \frac{n-1}{1-t} - (n+2) \Big) \ b_1''(0) $$
for all $b$ of the above type. With $b_1''(0) < 0$, we infer that $e_1$ is a local maximum of $\ln f$, $f$ and $A(.,t)$ if $\frac{n-1}{1-t} > n+2$, i.e. if $t > \frac 3 {n+2}$, and $e_1$ is a local minimum if $0 < t < \frac 3 {n+2}$. The neighborhood $\mathcal{U} \subset S^{n-1}$ of $e_1$ has to be chosen such that $a \in \mathcal{U}$ satisfies $a_1 > t > |a_j|$ for all $j=2 \etc n$. Since $e_1$ was the only critical point for $t > \frac 1 {\sqrt 2}$ with $a_1 > t > a_j$, this also gives another proof of \eqref{eq3.2}, i.e. the result of \cite{LT}.         \hfill $\Box$
\vspace{0,4cm}

An easy explicit example that for $0 < t = \frac 2 n < \frac 3 {n+2}$, $e_1$ is no longer a maximum of $A(.,t)$, is given by the vector
$\tilde{a} = (1 - \frac 2 n, \frac 2 n , \dots , \frac 2 n) \in S^{n-1} \subset \R^n$ where
$A(\tilde{a},\frac 2 n) = (\frac{n-2} n)^{n-2} > A(e_1,\frac 2 n) = (\frac{n-2} n)^{n-1}$. A continuous curve in $S^{n-1}$ joining $e_1$ and $\bar{a}$ shows that $e_1$ is not a local maximum for $A(.,\frac 2 n)$. \\

\remark The proof showed that for $\frac 1 {\sqrt 2} < t \le 1$ and for $0 < t < \frac 3 {n+2}$, $n \ge 6$, $e_1$ is the only critical point of $f( . ,t)$, a maximum in the first case and a minimum in the second. Therefore $e_1$ is the absolute maximum of $f(.,t)$ and the parallel section function $A(.,t)$ when
$t > \frac 1 {\sqrt 2}$, i.e. \eqref{eq3.2} holds, and $e_1$ is the absolute minimum of $f(.,t)$ when $0 < t < \frac 3 {n+2}$ and $a_1 > t > a_2 \etc a_n \ge 0$. However, this does not show that $e_1$ is an absolute minimum of $A(.,t)$ for $0 < t < \frac 3 {n+2}$ since formula \eqref{eq3.1} does not give the correct value of $A(a,t)$ for all $a \in S^{n-1}$ if $ t < \frac 1 {\sqrt 2}$; we have to assume the separation property $a_1 > t > a_j$, $j=2 \etc n$ for \eqref{eq3.1} to hold. Thus we can only claim ${e_1}$ to be locally extremal for $A(.,t)$. \\

\vspace{0,5cm}
To prove Theorem \ref{th3} on the maximal perimeter of non-central sections of the cross-polytope, we first show the following formula for the perimeter function

\begin{proposition}\label{prop3.1}
Let $a \in S^{n-1}$, $t>0$ with $a_1 >t >a_2 \etc a_n \ge 0$. Then
\begin{equation}\label{eq3.3}
P(a,t) = \frac{\sqrt n}{(n-2)!} \ \sum_{\e \in \{-1,+1\}^{n-1}} \sqrt{1 - \frac 1 n \pr{a}{\e}^2 } \ \prod_{j=2}^n \frac 1 {a_1 - \e_j a_j} \ (a_1 - t)^{n-2} \ ,
\end{equation}
where $\pr{a}{\e} = \sum_{j=1}^n \e_j a_j$ with $\e_1:=1$ and $\e_j \in \{-1,+1\}$.
\end{proposition}

\begin{proof}
There are $(2n-2)$ edges originating in $e_1$, from $e_1$ to $\e_j e_j$, $j=2 \etc n$, $\e_j \in \{-1,+1\}$. The hyperplane $H_t(a)$ intersects the boundary $\partial B(l_1^n)$ in $2^{n-1}$ faces $F_\e$, $F_\e := \text { convex hull of \ } ( e_1, \e_2 e_2 \etc \e_n e_n )$ with $\e = (\e_2 \etc \e_n) \in \{-1,+1\}^{n-1}$. We will determine $\vol_{n-2}(H_t(a) \cap F_\e)$, i.e. the volume of a hyperplane section of the $(n-1)$-simplex $F_\e$. The sum of all volumes over all
$\e = (\e_2 \etc \e_n) \in \{-1,+1\}^{n-1}$ will be $P(a,t)$. Let $P_\e$ be the pyramid - i.e. the irregular $(n-1)$-simplex - spanned by $e_1$ and $H_t(a) \cap F_\e$, i.e. $P_\e = \text{ convex hull of } ( e_1, v_2 \etc v_n )$ where $v_j = s_j e_1 + (1-s_j) \e_j e_j \in H_t(a) \cap span(e_1, \e_j e_j)$. The vectors
$v_j-e_1 =(1-s_j) (\e_j e_j - e_1)$ span $P_\e$. Since $t = \pr{a}{v_j} = s_j a_1 + (1-s_j) \e_j a_j$, $(1-s_j) = \frac{a_1-t}{a_1- \e_j a_j}$. We get as in the proof of Proposition \ref{prop2.1}
$$vol_{n-1} (P_\e) = \prod_{j=2}^n (1-s_j) \ \vol_{n-1}(span(\e_j e_j - e_1)_{j=2 \etc n}) = \frac{\sqrt n}{(n-1)!} \prod_{j=2}^n \frac{a_1-t}{a_1- \e_j a_j} \ , $$
since $span(\e_j e_j - e_1)_{j=2 \etc n}$ is a regular $(n-1)$-simplex. Let $c := \frac 1 n (\e_j)_{j=1}^n$ denote the centroid of $F_\e$, with $\e_1 = 1$. Define
$$b := \frac{a - \pr{a}{\e} c }{\norm{a - \pr{a}{\e} c}_2} = \Big (\frac{a_j - \frac 1 n \pr{a}{\e} \e_j }{\sqrt{1- \frac 1 n \pr{a}{\e}^2 }} \Big)_{j=1}^n \in \R^n \ . $$
Then $\norm{b}_2 = 1$ and $\pr{b}{c} = 0$. We find that
$$ \pr{a}{x} = t \text { \; if and only if \; } \pr{b}{x} = \frac{t- \frac 1 n \pr{a}{\e} }{\sqrt{1- \frac 1 n \pr{a}{\e}^2 }} =: s \ . $$
Therefore $\pr{a}{v_j} = t$ if and only if $\pr{b}{v_j} = s$, $j = 2 \etc n$. Moreover, $\pr{b}{e_1} = b_1$. Hence the pyramid $P_\e$ has height
$b_1 - s = \frac{a_1-t}{\sqrt{1- \frac 1 n \pr{a}{\e}^2 }}$ and therefore
\begin{align*}
vol_{n-2}(H_t(a) \cap F_\e) &= \frac{n-1}{b_1-s} \vol_{n-1}(P_\e) = \frac{\sqrt n}{(n-2)!} \frac 1 {b_1 - s} \ \prod_{j=2}^n \frac{a_1-t}{a_1- \e_j a_j} \\
&= \frac{\sqrt n}{(n-2)!} \ \sqrt{1- \frac 1 n \pr{a}{\e}^2 } \ \prod_{j=2}^n \frac 1 {a_1 - \e_j a_j} \ (a_1 - t)^{n-2} \ .
\end{align*}
Summing up all these volumes over $\e \in \{-1,+1\}^{n-1}$ yields the formula for $P(a,t)$.
\end{proof}

\vspace{0,5cm}
{\it Proof of Theorem \ref{th3}.} \\
Let $a \in S^{n-1}$, $a_1 > t > a_1 \etc a_n \ge 0$. Proposition \ref{prop3.1} and the Cauchy-Schwarz inequality yield
\begin{align*}
P(a,t) & = \frac{\sqrt n}{(n-2)!} \ \sum_{\e \in \{-1,+1 \}^{n-1}} \sqrt{1 - \frac 1 n \pr{a}{\e}^2 } \ \prod_{j=2}^n \frac 1 {a_1 - \e_j a_j} \ (a_1 - t)^{n-2} \\
& \le \frac{\sqrt n}{(n-2)!} \ \Big( \ \sum_{\e \in \{-1,+1 \}^{n-1}} ( 1 - \frac 1 n \pr{a}{\e}^2) \ \Big)^{1/2} \
\Big( \ \sum_{\e \in \{-1,+1\}^{n-1}} \prod_{j=2}^n \frac 1 {(a_1 - \e_j a_j)^2} \ \Big)^{1/2} \ (a_1 - t)^{n-2} \\
&= \frac{\sqrt n}{(n-2)!} \ \Big(2^{n-1} (1 - \frac 1 n) \Big)^{1/2} \Big( \ \prod_{j=2}^n \ [ \frac 1 {(a_1-a_j)^2}+ \frac 1 {(a_1+a_j)^2} ] \ \Big)^{1/2} \ (a_1 - t)^{n-2} \\
&= \frac{\sqrt {n-1}}{(n-2)!} \ 2^{(n-1)/2} \ \Big( \ \prod_{j=2}^n \ \frac{ 2(a_1^2+a_j^2) }{(a_1^2+a_j^2)^2} \ \Big)^{1/2} \ (a_1 -t)^{n-2} \\
&= \frac{\sqrt {n-1}}{(n-2)!} \ 2^{n-1} \ \prod_{j=2}^n \frac {\sqrt{1 + (a_j/a_1)^2}}{1- (a_j/a_1)^2} \ \frac{(a_1 -t)^{n-2}}{a_1^{n-1}} \ .
\end{align*}
The function $g(x) := \frac{\sqrt{1+x}}{1-x}$ is log-convex since $(\ln g)''(x) = \frac{1+6 x +x^2}{2(1-x)^4} >0$, with $\ln g(0) = 0$. Therefore \\
$\sum_{j=2}^n (\ln g)(x_j) \le (\ln g)(\sum_{j=2}^n x_j)$ \ for \  $x_j := (\frac {a_j}{a_1})^2$, \ $\sum_{j=2}^n x_j = \sum_{j=2}^n (\frac{a_j}{a_1})^2 = \frac{1-a_1^2}{a_1^2}$ , \\
$\ln \ \prod_{j=2}^n \frac {\sqrt{1 + (a_j/a_1)^2}}{1- (a_j/a_1)^2} = \sum_{j=2}^n (\ln g)(x_j) \le (\ln g)(\sum_{j=2}^n x_j) = \ln g(\frac{1-a_1^2}{a_1^2}) = \ln( \frac {a_1}{2 a_1^2 -1 } )$, \\
$\frac {\sqrt{1 + (a_j/a_1)^2}}{1- (a_j/a_1)^2} \le \frac{a_1}{2 a_1^2 -1 }$. We conclude that
$$P(a,t) \le \frac{\sqrt{n-1}}{(n-2)!} \ 2^{n-1} \ \frac 1 {2 a_1^2 -1} \ (1- \frac t {a_1})^{n-2} \ .$$
Let $h(a_1) := \frac 1 {2 a_1^2 -1} (1-\frac t {a_1})^{n-2}$. Then $h'(a_1) = \frac{(1- t/a_1)^{n-3}}{a_1^2 (2 a_1^2-1)^2} \ ( t [n(2 a_1^2-1)+2] - 4 a_1^3 )$. For $n \ge 4$ and $a_1 \ge t \ge \frac 1 {\sqrt 2}$, $h'(a_1) \ge 0$; namely for $n \ge 5$
$$t [n(2 a_1^2-1)+2] - 4 a_1^3 \ge \frac 1 {\sqrt 2} \ [5(2 a_1^2-1)+2] - 4 a_1^2 = a_1^2 \ [5 \sqrt 2 - 4 a_1] - \frac 3 2 \sqrt 2 \ge 0 \  $$
and for $n=4$
$$t [n(2 a_1^2-1)+2] - 4 a_1^3 \ge \frac 1 {\sqrt 2} \ [4(2 a_1^2-1)+2] - 4 a_1^2 = a_1^2 \ [4 \sqrt 2 - 4 a_1] -  \sqrt 2 := \psi(a_1) \ .  $$
The function $\psi$ is increasing in $(\frac 1 {\sqrt 2}, \frac{2 \sqrt 2} 3)$ and decreasing in $(\frac{2 \sqrt 2} 3,1)$ with $\psi(\frac 1 {\sqrt 2})=0$,
$\psi(\frac{2 \sqrt 2} 3) = \frac 5 {27} \sqrt 2$ and $\psi(1) = 3 \sqrt 2 -4 >0$, i.e. $\psi(a_1) \ge 0$. Therefore for all $n \ge 4$, $t \ge \frac 1 {\sqrt 2}$, $h$ is increasing in $a_1 \in (t,1)$, hence $h(a_1) \le h(1)$. This means that
$$P(a,t) \le \frac{\sqrt{n-1}}{(n-2)!} \ 2^{n-1} \ (1-t)^{n-2} = P(e_1,t) \ , $$
where $P(e_1,t)$ is calculated from \eqref{eq3.3}. For $n=3$ and $t \ge \frac 4 5$, the same is true since
$$t [n(2 a_1^2-1)+2] - 4 a_1^3 \ge \frac 4 5 \ [6 a_1^2 -1] - 4 a_1^2 = a_1^2 \ [\frac{24} 5 - 4 a_1] -  \frac 4 5 := \phi(a_1) \ ,  $$
with $\phi(a_1)$ decreasing in $(\frac 4 5 , 1)$, $\phi(\frac 4 5) = \frac {28}{125}$, $\phi(1)=0$. Therefore $h$ is increasing for $n=3$, too, if $t > \frac 4 5$.     \hfill $\Box$

\vspace{0,5cm}
{\it Proof of Corollary \ref{cor1}.}
Let $n \ge 6$, $ \frac 4 n < t \le \frac 1 {\sqrt 2}$ and \\
$\mathcal{U} := \{ a \in S^{n-1} \ | \ a_1 > \sqrt{1 - \frac 1 n}, \ t > |a_j| , j=2 \etc n \}$. Then $\mathcal{U}$ is an open neighborhood of of $e_1$ and the previous proof of Theorem \ref{th3} applies to all $a \in \mathcal{U}$, provided that $h'(a_1) \ge 0$, i.e. that $t[n(2 a_1^2 -1) +2] - 4 a_1^3 \ge 0$. Assuming $a_1 > \sqrt{1 - \frac 1 n }$, this is true since
$$t[n(2 a_1^2 -1) +2] - 4 a_1^3 \ge \frac 4 n \ [ n(1 - \frac 2 n) +2 ] = 4 \ge 4 a_1^2 \ . $$
Therefore $A(a,t) \le A(e_1,t)$ for all $a \in \mathcal{U}$, and $e_1$ is a local maximum of the perimeter function $P(.,t)$.   \hfill $\Box$

\vspace{0,5cm}
\remark The question whether the central perimeter function $P( . , 0)$ of the cross-polytope attains its maximal value also in $e_1$ - which the central section function $A( . , 0)$ does by Meyer and Pajor \cite{MP} - is open. However, we can prove a slightly weaker result, which is asymptotically optimal.

\begin{proposition}\label{prop3.2}
For any $a \in S^{n-1}$ we have
$$P(a,0) \le \sqrt{\frac n {n-1}} \ P(e_1,0) \ . $$
\end{proposition}

\begin{proof}
For $a \in S^{n-1}$, consider the central section of the cross-polytope by $H_0(a) = a^\perp$. We first give a formula for the part of the perimeter in the face $F_\e$ of $B(l_1^n)$, $F_\e = $ span $(\e_1 e_1\etc \ \e_n e_n)$, where $\e = (\e_1 \etc \e_n) \in \{-1,+1\}^n$, i.e. a formula for
$\vol_{n-2}(H_0(a) \cap F_\e)$. The face $F_\e$ is an $(n-1)$-simplex and by Dirksen \cite{D}, Theorem 1.1, we have
$$\vol_{n-2}((H_0(a) \cap F_\e) = \frac 1 {(n-2)!} \ \sqrt{n - \pr{a}{\e}^2 } \ \frac 1 {2 \pi} \int_\R \prod_{j=1}^n \ \frac 1 {1 + i \e_j a_j s} \ ds \ , $$
with $\pr{a}{\e} = \sum_{j=1}^n a_j \e_j$. Therefore
$$P(a,0) = \frac 1 {(n-2)!} \ \sum_{\e \in \{-1,+1\}^n} \sqrt{n - \pr{a}{\e}^2 } \ \frac 1 {2 \pi} \int_\R \prod_{j=1}^n \ \frac 1 {1 + i \e_j a_j s} \ ds \ . $$
Note that all integrals are non-negative since they are proportional to $(n-2)$-dimensional volumes. They may be zero, if $H_0(a)$ does not intersect $F_\e$. Estimating the weight by $\sqrt n$, we find
\begin{align*}
P(a,0) & \le \frac{\sqrt n}{(n-2)!} \ \frac 1 {2 \pi} \ \int_\R \sum_{\e \in \{-1,+1\}^n} \prod_{j=1}^n \ \frac 1 {1 + i \e_j a_j s} \ ds  \\
& = \frac{\sqrt n}{(n-2)!} \ \frac 1 {2 \pi} \ \int_\R  \ \prod_{j=1}^n \Big( \frac 1 {1 + i a_j s} + \frac 1 {1 - i a_j s} \Big) \ ds \\
& = \frac{\sqrt n}{(n-2)!} \ 2^{n-1}  \frac 1 {\pi} \ \int_\R \ \prod_{j=1}^n  \frac 1 {1 + a_j^2 s^2} \ ds  \ .
\end{align*}
By Meyer and Pajor \cite{MP}, Lemma II.6 and Proposition II.7, the integral is proportional to $A(a,0)$ and attains its maximum for $a=e_1$,
$$\frac 1 {\pi} \ \int_\R \ \prod_{j=1}^n  \frac 1 {1 + a_j^2 s^2} \ ds = \frac{(n-1)!}{2^{n-1}} A(a,0) \le \frac{(n-1)!}{2^{n-1}} A(e_1,0) = 1 \ . $$
Hence $P(a,0) \le \frac 1 {(n-2)!} \ \sqrt n \ 2^{n-1}$, while by the above formula for $P$, $P(e_1,0) = \frac 1 {(n-2)!} \ \sqrt{n-1} \ 2^{n-1}$, so that the claim $P(a,0) \le \sqrt{\frac n {n-1}} \ P(e_1,0)$ follows.
\end{proof}

\vspace{0,5cm}

\section{Non-central sections of the cube}

Let $Q^n := [-\frac 1 2, \frac 1 2]^n \subset \R^n$ be the cube of volume 1 and a distance $t$ be given with $\frac{\sqrt{n-1}} 2 < t \le \frac{\sqrt n} 2$. Denote the main diagonal by $a^{[n]} := \frac 1 {\sqrt n} (1 \etc 1) \in S^{n-1}$. As mentioned, Moody, Stone, Zach and Zvavitch \cite{MSZZ} showed that for all $a \in S^{n-1}$
\begin{equation}\label{eq4.1}
A(a,t) \le A(a^{[n]},t) = \frac{n^{n/2}}{(n-1)!} \ (\frac{\sqrt n} 2 - t)^{n-1} \ ,
\end{equation}
if $n \ge 3$. If $n=2$, $\frac 3 8 \sqrt 2 < t \le \frac 1 {\sqrt 2}$ is required. We extend this result to slightly smaller values of hyperplane distances $t$ to zero, namely to $\frac{n-k} 2 := \frac{n-2} {2 \sqrt n} < t \le \frac{\sqrt{n-1}} 2$, i.e. $k:= 4 - \frac 4 n$, in the sense that $a^{[n]}$ then is at least a {\it local} maximum. \\

We will use a different volume formula for $A(a,t)$ than the one used in \cite{MSZZ} to prove this result. The argument also yields a new proof of \eqref{eq4.1}. If $t > \frac {\sqrt{n-1}} 2$, $H_t(a) \cap Q^n$ separates one vertex of $Q^n$, say $e := \frac 1 2 (1 \etc 1) \in Q^n$, from all others, since if
$\pr{a}{e} \ge t$ and also $\pr{a}{e^{(i)}} \ge t$ for the vertex $e^{(i)} = \frac 1 2 (1 \etc 1,\underbrace{-1}_i, 1 \etc 1) \in Q^n$, we would get the contradiction \\
$\frac{\sqrt{n-1}} 2 < t \le \pr{a}{\frac{e+e^{(i)}} 2} \le \norm{\frac{e+e^{(i)}} 2}_2 = \frac{\sqrt{n-1}} 2$. Hence for all $i=1 \etc n$ \\
$\pr{a}{e} = \frac 1 2 \sum_{j=1}^n a_j \ge t > \pr{a}{e^{(i)}} = \frac 1 2 \sum_{j=1}^n a_j - a_i$, so that $a_i > 0$ follows. The following formula holds for all $a \in S^{n-1}$ and $t$ satisfying $\pr{a}{e} \ge t > \pr{a}{e^{(i)}}$ for all $i = 1 \etc n$. This means that the hyperplane $H_t(a)$ separates the vertex $e$ from all other vertices, in particular from $e^{(i)}$.

\begin{proposition}\label{prop4.1}
Let $K= Q^n$, $a \in S^{n-1}$ and $t>0$ be such that $\pr{a}{e} \ge t > \pr{a}{e^{(i)}}$ for all $i = 1 \etc n$. Then
\begin{equation}\label{eq4.2}
A(a,t) = \frac 1 {(n-1)!} \ \frac{(\frac 1 2 \sum_{j=1}^n a_j - t )^{n-1}}{\prod_{j=1}^n a_j} \ .
\end{equation}
\end{proposition}

\begin{proof}
Let $v_i$ denote the intersection of the hyperplane $H_t(a)$ with the edge from $e$ to $e^{(i)}$. Then $\{ x \in Q^n \ | \ \pr{a}{x} = t \}$ is the convex hull of $v_1 \etc v_n$. Let $P$ be the pyramid spanned by $e, v_1 \etc v_n$, i.e. the convex hull spanned by these vectors. Since $\pr{a}{e} = \frac 1 2 \sum_{j=1}^n a_j$ \ , \ $\pr{a}{v_j} = t$ \ , \ $\pr{a}{e-v_j} = \frac 1 2 \sum_{j=1}^n a_j - t$, hence $h := \frac 1 2 \sum_{j=1}^n a_j - t$ is the height of $P$. Therefore
$$A(a,t) = \vol_{n-1} ( \{ x \in Q^n \ | \ \pr{a}{x} = t \} ) = n \ \frac{\vol_n(P)}{\frac 1 2 \sum_{j=1}^n a_j - t} \ . $$
Let $v_j = s_j e + (1-s_j) e^{(j)}$. Then \\
$t = \pr{a}{v_j} = s_j (\frac 1 2 \sum_{i=1}^n a_i) + (1-s_j) (\frac 1 2 \sum_{i=1, i \ne j}^n a_i - \frac 1 2 a_j) =
\frac 1 2 \sum_{i=1}^n a_i - (1-s_j) a_j $. Hence $1-s_j = \frac{\frac 1 2 \sum_{i=1}^n a_i - t}{a_j}$ and $v_j - e = -(1-s_j)(e^{(j)}-e) = (1-s_j) e_j$ . Therefore
$$\vol_n(P) = \frac 1 {n!} \ \det ( (1-s_1) e_1 \etc (1-s_n) e_n ) = \frac 1 {n!} \prod_{j=1}^n (1-s_j) = \frac 1 {n!} \frac{(\frac 1 2 \sum_{i=1}^n a_i - t)^n}{\prod_{i=1}^n a_i} \ , $$
$$A(a,t) = \frac 1 {(n-1)!} \ \frac{(\frac 1 2 \sum_{i=1}^n a_i - t)^{n-1}}{\prod_{i=1}^n a_i} \ . $$
\end{proof}

\vspace{0,4cm}
This immediately implies that $A(a^{[n]},t) = \frac{n^{n/2}}{(n-1)!} \ ( \frac{\sqrt n} 2 - t )^{n-1}$ holds under the conditions of Proposition \ref{prop4.1}, which are satisfied, if $t > \frac{n-2} {\sqrt n}$. The centroid of the vertices $e^{(i)}$ adjacent to $e$ is
$c = \frac 1 n \sum_{i=1}^n e^{(i)} = \frac{n-2} n e$  and  \\
$\pr{a^{[n]}}{e} = \frac 1 2 \sqrt n > \pr{a^{[n]}}{e^{(i)}} = \pr{a^{[n]}}{c} = \frac 1 2 \frac{n-2}{\sqrt n} = \frac 1 2 \sqrt{n-k} \ , \ k:= 4 - \frac 4 n$.
Thus formula \eqref{eq4.2} will not extend to any $a \in S^{n-1}$ and $t$ with $t < \frac 1 2 \frac{n-2}{\sqrt n}$, since then $H_t(a)$ will no longer separate the vertex $e$ from all other vertices. We now prove Theorem \ref{th4} and Proposition \ref{prop4} on the (local) maxima of $A(. ,t)$. \\

\vspace{0,5cm}
{\it Proof of Theorem \ref{th4} and Proposition \ref{prop4}.} \\
(i) Suppose that $t>0$ satisfies $\sqrt{n-k} < 2 t \le \sqrt n$. In the case of Theorem \ref{th4}, i.e. \cite{MSZZ}, $k=1$, and in the case of Proposition \ref{prop4}, $k= 4 - \frac 4 n$, $\sqrt{n-k} = \frac{n-2}{\sqrt n}$. Let \\
$$\mathcal{U} :=\{ a \in S^{n-1} \ | \ \pr{a}{e} = \frac 1 2 \sum_{i=1}^n a_i > t > \frac 1 2 \sum_{i=1}^n a_i - a_j \ , \  \text{ \ for all \ } j=1 \etc n \} \ . $$
Recall that, if $k=1$, $t> \frac{\sqrt{n-1}} 2$, the requirement $t > \frac 1 2 \sum_{i=1}^n a_i - a_j$ for all $j$ is true automatically; for $k= 4 - \frac 4 n$ it is an additional condition needed for Proposition \ref{prop4.1} to hold. Consider
$$f(a,t):= \frac{(\sum_{i=1}^n a_i - 2 t)^{n-1}}{\prod_{i=1}^n a_i} \quad , \quad \sqrt{n-k} < 2 t \le \sum_{i=1}^n a_i \le \sqrt n $$
so that $A(a,t) = \frac 1 {2^{n-1} (n-1)!} f(a,t)$. We will show that $f(.,t)$ attains the absolute maximum in $\mathcal{U}$ in the vector $a^{[n]}$ in the case of Theorem \ref{th4} and that $a^{[n]}$ is a relative maximum in the case of Proposition \ref{prop4}. \\

Consider the Lagrange function of $\ln f$ relative to the constraint $\sum_{i=1}^n a_i^2 = 1$,
$$F(a,t) = (n-1) \ln (\sum_{i=1}^n a_i - 2 t) - \sum_{i=1}^n \ln a_i + \frac{\lambda}2 (\sum_{i=1}^n a_i^2 - 1) \ . $$
The equations for a critical point of $\ln f$ on $S^{n-1}$ are
\begin{equation}\label{eq4.3}
\frac{n-1}{\sum_{i=1}^n a_i - 2 t} - \frac 1 {a_j} + \lambda a_j = 0 \quad , \quad j= 1 \etc n \ .
\end{equation}
Let $S:= \sum_{i=1}^n a_i$. Multiplying the equations \eqref{eq4.3} by $a_j$ and adding them gives
$$(n-1) \frac S {S - 2 t} - n + \lambda = 0 \quad , \quad \lambda = - \frac{2 t n - S}{S - 2 t} \ . $$
Inserting this value of $\lambda$ into \eqref{eq4.3} and multiplying by $a_j$ yields the quadratic equation for $a_j$
$$a_j^2 - \frac{n-1}{2 t n - S} \ a_j + \frac{S - 2 t}{2 t n - S} = 0 \quad , \quad j= 1 \etc n \ , $$
with solutions
$$a_{j,\pm} = \frac 1 2 \frac {n-1}{2tn-S} \pm \sqrt{ (\frac 1 2 \frac {n-1}{2tn-S})^2 - \frac{S-2t}{2tn-S} } \ . $$
We need $a \in \mathcal{U}$, i.e. $S > 2 t > S - 2 a_j$ for all $j = 1 \etc n$. We claim that $a_{j,-}$ does not satisfy $2 t > S - 2 a_{j,-}$ if $n \ge 4$. Then all $a_j$ would have to be equal to the positive sign solution $a_{j,+}$ of the quadratic equation, whose coefficients are independent of $j$. The inequality
$2 t > S - 2 a_{j,-}$ is equivalent to
$$S - 2 a_{j,-} = S - \frac {n-1}{2tn-S} + \sqrt{ (\frac {n-1}{2tn-S})^2 - 4 \ \frac{S-2t}{2tn-S} } < 2 t \ , $$
i.e.
\begin{align*}
(\frac {n-1}{2tn-S})^2 & - 4 \ \frac{S-2t}{2tn-S} < \Big( \frac {n-1}{2tn-S} - (S - 2t) \Big)^2 \\
& = (\frac {n-1}{2tn-S})^2 - 2 \ \frac {n-1}{2tn-S} \  (S - 2t) + (S - 2t)^2 \ ,
\end{align*}
which means $2(n-3) < (S-2t)(2tn-S)$ or
\begin{equation}\label{eq4.4}
S^2 - 2 (n+1) t S + 4 t n^2 + 2(n-3) < 0 \ .
\end{equation}
The corresponding quadratic {\it equation} has the solutions
$$S_{\pm} = (n+1) t \pm \sqrt{((n-1) t)^2 - 2(n-3)} \ . $$ Thus, if $S$ satisfies the quadratic {\it inequality} \eqref{eq4.4}, we have $S_- < S < S_+$, in particular $S> S_-$. Now $S_- = (n+1) t - \sqrt{((n-1) t)^2 - 2(n-3)}$ is increasing as a function of $t$, the $t$-derivative being
$n+1 - \frac{(n-1)^2 t}{\sqrt{((n-1) t)^2 - 2(n-3)}}$. This is positive if equivalently $4 n (n-1)^2 t^2 > 2 (n-3)(n+1)^2$ and, with $2 t \ge \frac{n-2}{\sqrt n}$, i.e. $4 n t^2 \ge (n-2)^2$, the positivity follows from $(n-1)^2 (n-2)^2 \ge 2 (n-3)(n+1)^2$, which is true for any $n \ge 3$. Replacing $t$ by the smaller $\frac{n-2}{2 \sqrt n}$, we find for any $n \ge 4$
\begin{align*}
S_- \ge & \frac{(n+1)(n-2)}{2 \sqrt n} - \sqrt{ \frac{(n+1)^2(n-2)^2}{4 n} - 2 (n-3) } \\
& = \frac{(n+1)(n-2)}{2 \sqrt n} - \sqrt{ \frac{(n^2-3 n - 2)^2}{ 4 n } } = \sqrt n \ .
\end{align*}
Note that the last equality holds only for $n \ge 4$ and not for $n=3$ since then $n^2-3n-2$ is negative. Hence for any $n \ge 4$, $S > S_- \ge \sqrt n$. However, this contradicts $ S = \sum_{i=1}^n a_i \le \sqrt n$. Therefore $a_{j,-}$ is not a possible solution of \eqref{eq4.3} and all coordinates of a critical point $a$ have to be equal (to $a_{j,+}$). Since $\sum_{i=1}^n a_i^2 =1$, $a_1 = \dots = a_n = \frac 1 {\sqrt n}$, and $a = a^{[n]}$ for $n \ge 4$ is the only critical point of $f(.,t)$ in $\mathcal{U}$. Clearly $a^{[n]}$ satisfies \eqref{eq4.3} with appropriate $\lambda$. We now show that $a^{[n]}$ is a relative maximum of $f( ., t)$. Then $a^{[n]}$ is the absolute maximum of $f(.,t)$ and $A(.,t)$ in the case of Theorem \ref{th4} ($k=1$) and a relative maximum in the case of Proposition \ref{prop4} ($k=4-\frac 4 n$). \\

\vspace{0,3cm}
(ii) Take a smooth non-constant $C^2$-function $b : [-\e_0,\e_0] \to \mathcal{U}$ with $b(0) = a^{[n]}$. We will show that $\ln f( . ,t) \circ b$ has a local maximum in zero for all such functions $b$, with $t$ fixed. Then $\ln f$ and $f$ will have a local maximum in $a^{[n]}$. Since $\sum_{i=1}^n b_i^2 = 1$,
$$\sum_{i=1}^n b_i b_i' = 0 \quad , \quad \sum_{i=1}^n(b_i'^2 + b_i b_i'') = 0 \ . $$
Clearly $(\ln f \circ b)'(0) = 0$, since $a^{[n]}$ is a critical point of $f$. Calculation shows
$$(\ln f \circ b)'' = - \frac{n-1}{(\sum_{i=1}^n b_i - 2t)^2} \ (\sum_{i=1}^n b_i')^2 + \frac{n-1}{\sum_{i=1}^n b_i - 2t} \ \sum_{i=1}^n b_i''
- \sum_{i=1}^n \frac{b_i''}{b_i} + \sum_{i=1}^n (\frac{b_i'}{b_i})^2 \ . $$
With $b(0) = \frac 1 {\sqrt n} (1 \etc 1)$, $\sum_{i=1}^n b_i'(0) = 0$ and $\sum_{i=1}^n b_i''(0) = - \sqrt n \sum_{i=1}^n b_i'(0)^2 < 0 $, we get
\begin{align*}
(\ln f \circ b)''(0) & = \frac{n-1}{\sqrt n - 2 t} \sum_{i=1}^n b_i''(0) - \sqrt n \sum_{i=1}^n b_i''(0) + n \sum_{i=1}^n b_i'(0)^2 \\
& = \Big( \frac{n-1}{\sqrt n - 2 t} - 2 \sqrt n \Big) \ \sum_{i=1}^n b_i ''(0) < 0 \ ,
\end{align*}
if $\frac{n-1}{\sqrt n - 2 t} > 2 \sqrt n$, which means $ 2 t > \frac{n+1}{2 \sqrt n} $. However, $2 t > \frac{n-2}{\sqrt n} \ge \frac{n+1}{2 \sqrt n}$ for all $n \ge 5$. Therefore $a^{[n]}$ is a relative maximum of $f(.,t)$ and $A(.,t)$ for all $t$ considered in Theorem \ref{th4} and Proposition \ref{prop4}, if $n \ge 5$. For $n=4$, we have a relative maximum if $\frac 5 8 < t \le 1$, and also for $n=3$, if $\frac 1 {\sqrt 3} < t \le \frac {\sqrt 3} 2$. This covers also the cases $n=3, 4$ in Theorem \ref{th4} with $\frac{\sqrt{n-1}} 2 < t \le \frac{\sqrt n} 2$. However, for $n=4$ and $\frac 1 2 < t < \frac 5 8$, $f(.,t)$ and $A(.,t)$ have a relative minimum in $a^{[n]}$. For $n=3$, we also have a relative minimum if $\frac 1 {2 \sqrt 3} < t < \frac 1 {\sqrt 3}$. \\

\vspace{0,3cm}
(iii) For $n=3$ in Theorem \ref{th4} ($k=1$), there were possibly two choices for the coordinates $a_j$ of a critical point of $f(.,t)$ in part (i). Assume that $a=(a_1,a_1,a_2)$ with $a_2=\sqrt{1-2 a_1^2}$, $a_1 \in [0, \frac 1 {\sqrt 2}]$. Then
$$f(a,t) = \frac{(2 a_1 + \sqrt{1-2 a_1^2} - 2 t)^2}{a_1^2 \sqrt{1-2a_1^2}} \text{ \ with \ } 2 a_1 + \sqrt{1-2 a_1^2} > 2 t > \sqrt 2 \ . $$
We have $\frac{\partial f}{\partial {a_1}}(a_1,t) = \frac{2 (2 a_1 + \sqrt{1-2 a_1^2} - 2 t)}{(a_1 \sqrt{1-2 a_1^2})^3} \Big[2 a_1^3 - (1-a_1^2)\sqrt{1-2a_1^2} + t (1-3 a_1^2) \Big]$.
Let $\phi(a_1,t):= 2 a_1^3 - (1-a_1^2)\sqrt{1-2a_1^2} + t (1-3 a_1^2)$. We claim that for any $t \in [\sqrt 2, \sqrt 3]$, $\phi(a_1,t) > 0$ for $a_1 \in [0, \frac 1 {\sqrt 3})$ and $\phi(a_1,t) < 0$ for $a_1 \in (\frac 1 {\sqrt 3},\frac 1 {\sqrt 2})$. Define $g: [0,\frac 1 {\sqrt 2}] \to \R$ by
$g(a_1) := \frac{2a_1^3-(1-a_1^2) \sqrt{1-2 a_1^2}}{3 a_1^2 -1}$. Then $g$ is increasing and has a removable singularity at $\frac 1 {\sqrt 3}$, with $g(0)=1 < g(\frac 1 {\sqrt 3}) = \sqrt{\frac 4 3} < g(\frac 1 {\sqrt 2}) = \sqrt 2$. Hence for all $a_1 \in [0, \frac 1 {\sqrt 3})$,
$$(1-a_1^2) \sqrt{1-2 a_1^2} - 2 a_1^3 = (1-3 a_1^2) g(a_1) \le (1-3 a_1^2) g(\frac 1 {\sqrt 3}) = \sqrt{\frac 4 3} (1- 3 a_1^2) < t (1- 3 a_1^2) \ , $$
thus $g|_{[0, \frac 1 {\sqrt 3})} > 0$. For $a_1 \in (\frac 1 {\sqrt 3},\frac 1 {\sqrt 2})$,
$$2a_1^3-(1-a_1^2) \sqrt{1-2 a_1^2} = (3 a_1^2 -1) g(a_1) \le (3 a_1^2 -1) g(\frac 1 {\sqrt 2}) = \sqrt 2 (3 a_1^2 -1) \le t (3 a_1^2 -1) \ , $$
hence $g|_{(\frac 1 {\sqrt 3},\frac 1 {\sqrt 2})} < 0$. This implies that $f(.,t)|_{[0, \frac 1 {\sqrt 3})}$ is increasing and $f(.,t)|_{(\frac 1 {\sqrt 3},\frac 1 {\sqrt 2})}$ is decreasing for all $t \in [\sqrt 2, \sqrt 3]$. Hence $f(.,t)$ has an absolute maximum in $a_1 = \frac 1 {\sqrt 3}$, i.e. $a = a^{[3]}$. This proves Theorem \ref{th4} also in the case $n=3$.  \hfill $\Box$

\vspace{0,5cm}
\remark We quickly sketch the idea of another possible proof of Theorem \ref{th4} and Proposition \ref{prop4}. Take $a \in \mathcal{U}$ with $a \ne a^{[n]}$, put $a_1 = \max_j a_j$ and suppose that for a suitable coordinate, say $a_2$, we have $a_1 \ne a_2$. Define a new vector $\tilde{a}$ by
$a_1' = a_2' = \frac{a_1+a_2} 2$, $a_j' = a_j$ for $j \ge 3$ and $\tilde{a} = \frac{a'}{\norm{a'}_2}$ , $\norm{a'}_2 = \sqrt{1 - \frac 1 2 (a_1-a_2)^2} =: W$.
We claim that $A(a,t) < A(\tilde{a},t)$. If this is true, averaging will increase the volume $A(.,t)$, and the maximum will occur for the fully averaged main diagonal vector $a^{[n]}$. The claim is equivalent to
$$1 < \Big( 1 + \frac{2 t}{\sum_{j=1}^n a_j-2 t} \ (1 - W) \Big)^{n-1} \ W \ ( 1 - (\frac{a_1-a_2}{a_1+a_2})^2 ) \ . $$
With $x := \frac 1 2 (a_1 - a_2)^2$,
$\frac{2 t}{\sum_{j=1}^n a_j - 2 t} \ge A:= \left\{
\begin{array}{c@{\quad}l} (n-1) (\sqrt{\frac n {n-1}} +1 )  \quad  k=1 \\
\frac n 2 - 1 \quad \quad \quad \quad \quad k= 4 - \frac 4 n
\end{array} \right\} $
and $M:= \frac 2 {(a_1+a_2)^2}$, it then suffices to verify
$$1 < (1 + A \frac x 2)^{n-1} \ \sqrt{1-x} \ ( 1 - M x ) \ . $$
To do so, show for $f(x) := (1+A \frac x 2)^{n-1} \ \sqrt{1-x}$ that $f(x) \ge 1 + N x$ where \\
$N := \left\{
\begin{array}{c@{\quad}l} (n-1)^2  \quad \quad  k=1 \\
\frac{n(n-3)} 4 \quad  k= 4 - \frac 4 n
\end{array} \right\} $.
Then it remains to prove $N - M > N M \frac 1 2 (a_1-a_2)^2$ which is equivalent to $a_1 a_2 > \frac 1 {2 N}$. If $k=1$, it is possible to find $a_2 \ne a_1$ with this property, using the assumption $\sqrt{n-1}< 2 t \le \sum_{j=1}^n a_j \le \sqrt n$. In the second case $k=4 - \frac 4 n$, the neighborhood $\mathcal{U}$ of $a^{[n]}$ has to be restricted to prove $a_1 a_2 > \frac 1 {2 N}$. The condition means that $a_2$ has to be chosen not too small compared to $a_1$. \\
A similar proof can be given in the case of the cross-polytope results, but then it is better to use geometric mean averaging, i.e. $a_1' = a_2' = \sqrt{a_1 a_2}$.     \hfill $\Box$

\vspace{0,5cm}
To prove Theorem \ref{th5}, we first show a formula for the perimeter.

\begin{proposition}\label{prop4.2}
Let $K=Q^n$, $a \in S^{n-1}$ and $t>0$ be such that $\pr{a}{e} > t > \pr{a}{e^{(i)}}$ for all $i = 1 \etc n$. Then
$$P(a,t) = \frac 1 {(n-2)!} \ \sum_{k=1}^n a_k \sqrt{1-a_k^2} \  \frac{( \frac 1 2 \sum_{i=1}^n a_i - t )^{n-2}}{\prod_{i=1}^n a_i} \ . $$
\end{proposition}
The assumption is satisfied, in particular, when $\pr{a}{e} > t > \frac{\sqrt{n-1}} 2$. \\

\begin{proof}
The intersection of the hyperplane $H_t(a)$ with the boundary of the cube $Q^n$ consists of the intersection with $n$ faces, $H_t(a) \cap \partial Q^n = \cup_{i=1}^n F_i$, \\
$F_i = \{ x \in Q^n \ | \ x_i = \frac 1 2 \}$. Note that $H_t(a) \cap  G_i$, $G_i = \{ x \in Q^n \ | \ x_i = - \frac 1 2 \}$ is empty since e.g. if $x_1 = - \frac 1 2$, $x \in H_t(a) \cap G_1$, we would get the contradiction
$t = \pr{a}{x} \le - \frac 1 2 a_1 + \sum_{j=2}^n a_j = \frac 1 2 \sum_{j=1}^n a_j -a_1 = \pr{a}{e^{(i)}} < t $. The boundary face $F_1$ is centered at
$c:=\frac 1 2 (1,0 \etc 0)$ and $x=(\frac 1 2, \tilde{x}) \in H_t(a) \cap F_1$ satisfies $\pr{\tilde{a}}{\tilde{x}} = \pr{a}{x} - \frac 1 2 a_1 = \frac{2 t -a_1} 2$, with $\tilde{a}=(a_2 \etc a_n)$. Normalize $\tilde{a}$ by $\tilde{b} := \frac{\tilde{a}}{\sqrt{1-a_1^2}}$. Then
$\pr{\tilde{b}}{\tilde{x}} = \frac{ 2 t - a_1}{2 \sqrt{1-a_1^2}} =: s$ is the distance from the center $c$ of the face $F_1$ to $H_t(a) \cap F_1$. Using the assumption on $t$, we have for $i = 1$
$$\pr{\tilde{b}}{\tilde{e}} = \frac{\frac 1 2 \sum_{j=2}^n a_j } {\sqrt{1-a_1^2}} > \frac{ t - \frac 1 2 a_1}{\sqrt{1-a_1^2}} =s > \frac{\frac 1 2 \sum_{j=2}^n a_j - a_1} {\sqrt{1-a_1^2}} = \pr{\tilde{b}}{\tilde{e}^{(1)}} \ . $$
Therefore the assumption of Proposition \ref{prop4.1} is satisfied for the $(n-1)$-face $F_1$ and distance $s$. Hence
\begin{align*}
\vol_{n-2}(H_t(a) \cap F_1) &= \frac 1 {(n-2)!} \ \frac{(\frac 1 2 \sum_{j=2}^n \tilde{b}_j - s)^{n-2}}{\prod_{j=2}^n \tilde{b}_j} \\
&= \frac 1 {(n-2)!} \ \sqrt{1-a_1^2} \ \frac{(\frac 1 2 \sum_{j=2}^n a_j - (t- \frac 1 2 a_1))^{n-2}}{\prod_{j=2}^n a_j} \\
&= \frac 1 {(n-2)!} \ a_1 \sqrt{1-a_1^2} \ \frac{(\frac 1 2 \sum_{j=1}^n a_j - t)^{n-2}}{\prod_{j=1}^n a_j} \ .
\end{align*}
Summing the intersection volumes over all faces, we get
$$P(a,t) = \frac 1 {(n-2)!} \ \sum_{k=1}^n a_k \sqrt{1-a_k^2} \ \frac{(\frac 1 2 \sum_{j=1}^n a_j - t)^{n-2}}{\prod_{j=1}^n a_j} \ . $$
Let us mention that, if $t > \frac 1 2 \sqrt{n-1}$, then $s > \frac 1 2 \sqrt{n-2}$, since
$0 \le (\sqrt{n-1} \ a_1 -1 )^2$ implies $\sqrt{n-1} - a_1 \ge \sqrt{n-2} \sqrt{1-a_1^2}$ which gives $\frac 1 2 \sqrt{n-2} \le \frac{\sqrt{n-1} - a_1}{2 \sqrt{1-a_1^2}}$ and the latter is $< s$ if $t > \frac 1 2 \sqrt{n-1}$.
\end{proof}

\vspace{0,5cm}
{\it Proof of Theorem \ref{th5} and Proposition \ref{prop5}.} \\
(1) Let $n \ge 4$, $\sqrt{n-1} < 2 t \le \sqrt n$ (case $k=1$, Theorem \ref{th5}) or $n \ge 5$, $\sqrt{n-k} = \frac{n-2}{\sqrt n} < 2 t \le \sqrt{n-1}$ (case $k=4-\frac 4 n$, Proposition \ref{prop5}). We claim that $a = a^{[n]}$ maximizes the function $f( . , t)$,
$$f(a,t) := \frac{ (\sum_{i=1}^n a_i - 2 t)^{n-2}}{\prod_{i=1}^n a_i} \ , $$
on $\mathcal{U}$, where $\mathcal{U}$ is as in the proof of Theorem \ref{th4} and Proposition \ref{prop4}. The function $f(a, t)$ which differs from the formula for the function $f(a,t)$ there by the exponent, $(n-1)$ being replaced by $(n-2)$, and thus also by the difference of degrees of homogeneity of the numerator and the denominator. We proceed as in the proof of Theorem \ref{th4} and Proposition \ref{prop4}. The equations for a critical point of $f(.,t)$ relative to $S^{n-1}$ are now
$$\frac{n-2}{\sum_{i=1}^n a_i - 2 t} - \frac 1 {a_j} + \lambda a_j = 0 \quad , \quad j= 1 \dots n \ , $$
which yield with $S := \sum_{i=1}^n a_i$ that $\lambda = - \frac{2 t n - 2 S}{S - 2 t}$, giving the quadratic equation for the $a_j$'s
\begin{equation}\label{eq4.5}
a_j^2 - \frac{n-2} {2 n t - 2 S} \ a_j + \frac{S - 2 t}{2nt - 2S} = 0 \ ,
\end{equation}
$$a_{j , \pm} = \frac 1 4 \frac{n-2}{n t - S} \pm \sqrt{ (\frac 1 4 \frac{n-2}{n t - S})^2 - \frac{S-2t}{2(n t - S)} } \ . $$
Again we show that $a_{j,-}$ is not a possible solution since the requirement $2 t > S - 2 a_{j,-}$ will be violated if $n \ge 5$:
$$S - 2 a_{j,-} = S - \frac 1 2 \frac {n-2}{nt-S} + \sqrt{ (\frac 1 2 \frac {n-2}{nt-S})^2 - 2 \ \frac{S-2t}{nt-S} } < 2 t \ , $$
i.e.
\begin{align*}
(\frac 1 2 \frac {n-2}{nt-S})^2 & - 2 \ \frac{S-2t}{nt-S} < \Big( \frac 1 2 \frac {n-2}{nt-S} - (S - 2t) \Big)^2 \\
& = (\frac 1 2 \frac {n-2}{nt-S})^2 -  \frac {n-2}{nt-S} \  (S - 2t) + (S - 2t)^2 \ ,
\end{align*}
which is equivalent to $n-4 < (S-2 t) ( n t -S)$, i.e.
$$S^2 - (n+2) t S + 2 n t^2 + (n-4) < 0 \ . $$
The solutions of the corresponding quadratic equation are \\
$S_{\pm} = \frac{n+2}2  t \pm \sqrt{ (\frac{n-2} 2 t)^2 - (n-4) } $ and $S_- < S < S_+$ follows. Again $S_-$ is an increasing function of $t$ and using $t \ge \frac{n-2}{2 \sqrt n}$, we find for all $n \ge 5$
$$S_- \ge \frac{n^2-4}{4 \sqrt n} - \sqrt{ \frac{(n-2)^4}{16 n} - (n-4) } = \frac{n^2-4}{4 \sqrt n} - \sqrt{ \frac{(n^2-4 n - 4)^2}{16 n} } = \sqrt n \ . $$
The last equality does not hold for $n=4$ since then $n^2-4n-4$ is negative. Hence for $n \ge 5$, the solution $a_{j,-}$ of \eqref{eq4.5} yields no allowed solution $a\in \mathcal{U}$ since $S\sqrt n \ge S > S_- \ge \sqrt n$ would give a contradiction. Therefore the only critical point of $f(.,t)$ in $\mathcal{U}$ is $a^{[n]}$. We now show that it is a relative maximum if $n \ge 6$. \\

\vspace{0,3cm}
(ii) Proceeding as in part (ii) of the proof of Theorem \ref{th4} and Proposition \ref{prop4} we find similarly
$$(\ln f \circ b)''(0) = \Big( \frac{n-2}{\sqrt n - 2 t} - 2 \sqrt n \Big) \ \sum_{i=1}^n b_i''(0) < 0 \ , $$
if $\frac{n-2}{\sqrt n - 2 t} > 2 \sqrt n$, i.e. $2 t > \frac{n+2}{2 \sqrt n}$. We know $2 t > \frac{n-2}{\sqrt n}$ and this is $\ge \frac{n+2}{2 \sqrt n}$ for all $n \ge 6$. For $n=5$, it is true as long as $\frac 7 {2 \sqrt 5} < 2 t \le \sqrt 5$, for $n=4$, if $\frac 3 2 < 2 t \le 2$ and for $n=3$, if $\frac 5 {2 \sqrt 3} < 2 t \le \sqrt 3$: these conditions suffice to show that $a^{[n]}$ is a relative maximum of $A(.,t)$. This covers the range of $t$-values in Theorem \ref{th5} and in Proposition \ref{prop5} completely. \\

\vspace{0,3cm}
(iii) For $n=4$ we have possibly two solutions $a_{j,\pm}$ of equation \eqref{eq4.5} for a critical point of $f(.,t)$. Then either $a=(a_1,a_1,a_1,a_2)$ with $a_2=\sqrt{1-3a_1^2}$ or $a=(a_1,a_1,a_2,a_2)$ with $a_2=\sqrt{\frac{1-a_1^2} 2}$ and functions to maximize
$$ f_1(a_1,t) = \frac{(3 a_1 + \sqrt{1-3 a_1^2} - 2 t)^2}{a_1^3 \sqrt{1-3a_1^2}} \ , \ f_2(a_1,t) = 2 \frac{(2 a_1 + \sqrt{2 - 4 a_1^2} - 2 t)^2}{a_1^2 (1- 2 a_1^2)} \ . $$
For any $t \in [\sqrt 3 , 2]$, the investigation is similar to the corresponding function $f(.,t)$ in part (iii) of the proof of Theorem \ref{th4}, where $n$ was $3$. The functions $\phi_1$ and $\phi_2$ determining the sign of the derivative of $f_1$ and $f_2$, respectively, are
$$\phi_1(a_1,t) = 6 a_1^3 -a_1 - (1- 2 a_1^2) \ \sqrt{1-3 a_1^2} + (1-4 a_1^2 -1) t \ , $$
$$\phi_2(a_1,t) = 4 a_1^3 \ \sqrt{2-4 a_1^2} + 8 a_1^2 -8 a_1^4 -2 + \sqrt{2-4 a_1^2} \ (1-4 a_1^2) t \ . $$
The maximum of $f_1(.,t)$ and of $f_2(.,t)$ is attained for $a_1=\frac 1 2$ in both cases, yielding the unique absolute maximum
$a^{[4]}= \frac 1 2 (1,1,1,1)$ for both $f_i(.,t)$ and all $t \in [\sqrt 3 , 2]$. We do not give the detailed calculations since they are analogous to part (iii) of the proof of Theorem \ref{th4}. We only mention that in the case of $f_2$, only $t \ge \frac 3 2$ is needed, while for maximizing $f_1$, $t \ge \sqrt 3$ is required. \\

In dimension $n=3$, the result of Theorem \ref{th5} is true for $0.7248 \simeq \frac{\sqrt 3 + 8 \sqrt 2} {18} =: t_0 < t \le \frac{\sqrt 3} 2$:  For any fixed $t$ and $a=(a_1,a_1,a_2)$ with $a_2 = \sqrt{1-2 a_1^2}$, \\
$f(a,t)= \frac{ 2 a_1 + \sqrt{1-2 a_1^2} - 2t }{a_1^2 \sqrt{1-2a_1^2}}$ may again be analyzed by differential calculus. The result is that for
$0.7251 \simeq \frac 1 6 (\sqrt{18+10 \sqrt 3} - \sqrt{6 - 2 \sqrt 3}) =: t_1 < t \le \frac{\sqrt 3} 2$, the only critical point of $f( . , t)$ is $a^{[3]}$, whereas for $t < t_1$ there is another critical point $\bar{a}$ with $a_1 \ne a_2$. For $t_0 < t < t_1$,
$f(\bar{a},t) < f(a^{[3]},t)$ is still true, whereas for $t<t_0$, $\bar{a}$ gives the maximum of $f( . , t)$. \\

\vspace{0,3cm}
(iv) Using Proposition \ref{prop4.2} and parts (i) to (iii), we have for all $a \in \mathcal{U}$ and $n \ge 4$ that
$$P(a,t) \le \frac 1 {2^{n-2} (n-2)!} \ \sum_{j=1}^n a_j \sqrt{1-a_j^2} \ f(a^{[n]},t) \ . $$
The function $g(x) := \sqrt{x (1-x)}$ is concave in $[0,1]$. $g''(x) = - \frac 1 4 (x (1-x))^{-3/2} < 0$. Therefore
$$\frac 1 n \sum_{j=1}^n g(a_j^2) \le g(\frac 1 n \sum_{j=1}^n a_j^2) = g(\frac 1 n) = \frac{\sqrt{n-1}} n \ . $$
Hence
$$P(a,t) \le \frac 1 {2^{n-2} (n-2)!} \ \sqrt{n-1} \ f(a^{[n]},t) = \frac{\sqrt{n-1}}{(n-2)!} \ n^{n/2} (\frac{\sqrt n} 2 - t)^{n-2} = P(a^{[n]},t) \ . $$

For $n=3$, part (iii) implies that the perimeter function $P(. , t)$ attains its absolute maximum in $a^{[3]}$ for $t_0 < t \le \frac{\sqrt 3} 2 $. The numerical evidence is that $P( . , t)$ is also maximized by $a^{[3]}$ for $\frac{\sqrt 2} 2 < t \le t_0$.      \\

A lengthy calculation shows that
$$\frac {\partial^2 P}{\partial a_1^2}(\frac 1 {\sqrt 3}, \frac 1 {\sqrt 3},t) =\frac {\partial^2 P}{\partial a_2^2}(\frac 1 {\sqrt 3}, \frac 1 {\sqrt 3},t) = \frac 9 {2 \sqrt 2} (11 - 10 \sqrt 3 t) \ , $$
$$\frac {\partial^2 P}{\partial a_1 \partial a_2}(\frac 1 {\sqrt 3}, \frac 1 {\sqrt 3},t) = \frac 9 {4 \sqrt 2} (11 - 10 \sqrt 3 t) \ . $$
This implies that $a^{[3]}$ is a local {\it maximum} of the perimeter function $P( . , t)$ if $0.635 \simeq \frac {11}{30} \sqrt 3 < t < t_0$ and a local {\it minimum} if $0.577 \simeq \frac 1 {\sqrt 3} < t < \frac {11}{30} \sqrt 3 \simeq 0.635$. Hence the statement of Proposition \ref{prop5} is only partially true in dimension $n=3$.                 \hfill $\Box$

\vspace{0,5cm}
In the case of the cube, there is also a general formula for the perimeter at hyperplane distance $t$ from the origin, for any $t \in \R$, although it is less explicit than the one in Proposition \ref{prop4.2}.

\begin{proposition}\label{prop4.3}
Let $n \ge 3$, $K = Q^n$, $a \in S^{n-1} \subset \R^n$ and $t \in \R$. Then \\
(a) $A(a,t) = \frac 2 {\pi} \ \int_0^\infty \prod_{j=1}^n \frac{\sin(a_j s)}{a_j s} \cos(2 t s) \ ds \ , $ \\
(b) $P(a,t) = 2 \sum_{k=1}^n \sqrt{1-a_k^2} \ \frac 2 {\pi} \ \int_0^\infty \prod_{j=1,j \ne k}^n \frac{\sin(a_j s)}{a_j s} \ \cos(a_k s) \ \cos(2 t s) \ ds \ . $
\end{proposition}

\begin{proof}
Formula (a) for the parallel section function $A$ is well-known, cf. Ball \cite{B} or K\"onig and Koldobsky \cite{KK1}. In \cite{KK1} the distance to zero is $\frac t 2$, so $t$ there has to be replaced by $2t$ here. \\
To prove (b), consider the two faces $F_\pm :=\{ x \in Q^n \ | \ x_1 = \pm 1 \}$ of the cube and their hyperplane sections $H_t(a) \cap F_\pm$. The distance of the section $H_t(a) \cap F_\pm$ to the midpoint of $F_\pm$ is $t_\pm = \frac 1 2 \frac{2 t \mp a_1}{\sqrt{1-a_1^2}}$, since for
$x=(\pm \frac 1 2, \tilde{x}) \in F_\pm$ and $a=(a_1, \tilde{a})$ with $\norm{\tilde{a}}_2 = \sqrt{1-a_1^2}$, we have \
$ t = \pr{a}{x} = \pm \frac{a_1} 2 + \pr{\tilde{a}}{\tilde{x}}$ , $\pr{ \frac{\tilde{a}} {\norm{\tilde{a}}_2} }{\tilde{x}} = \frac{2 t \mp a_1}
{2 \norm{\tilde{a}}_2} $. \\
Therefore the part of the perimeter $P(a,t)$ on both faces $F_\pm$ is given by
\begin{align*}
\frac 2 {\pi} & \int_0^\infty \prod_{j=2}^n \frac{\sin(\frac{a_j}{\norm{\tilde{a}}_2} s)}{\frac{a_j}{\norm{\tilde{a}}_2} s} \
[\ \cos(\frac{2 t \mp a_j}{\norm{\tilde{a}}_2} s) + \cos(\frac{2 t \pm a_j}{\norm{\tilde{a}}_2} s) \ ] \ ds \\
& = 2 \sqrt{1-a_1^2} \ \frac 2 {\pi} \ \int_0^\infty \prod_{j=2}^n \frac{\sin(a_j r)}{a_j r} \ \cos(a_1 r) \ \cos(2 t r) \ dr \ .
\end{align*}
Summing up the contributions of all faces, we get
$$P(a,t) = 2 \sum_{k=1} \sqrt{1-a_k^2} \ \frac 2 {\pi} \ \int_0^\infty \prod_{j=1,j \ne k}^n \frac{\sin(a_j s)}{a_j s} \ \cos(a_k s) \ \cos(2 t s) \ ds \ . $$
\end{proof}

Note that some of the integrals in (b) may be zero if the hyperplane $H_t(a)$ does not intersect the corresponding faces. \\

In the case of central sections of the cube $Q^n$ we have the well-known inequalities foe all $a \in S^{n-1}$
$$ 1 \le A(e_1,0) \le A(a,0) \le A(a^{[2]},0) = \sqrt 2 \ , \ a^{[2]} = \frac 1 {\sqrt 2} (1,1,0 \etc 0) \ , $$
the lower bound having been shown by Hadwiger \cite{H} and the upper bound by Ball \cite{B}. The perimeter formula of Proposition \ref{prop4.3} was used for $t=0$ by K\"onig and Koldobsky \cite{KK} to prove the analogue of Ball's result for the perimeter, that
$$ P(a,0) \le P(a^{[2]},0) $$
holds for all $a \in S^{n-1}$. The case $n=3$ of this estimate was first shown by Pelczy\'nski (private communication, 2005). \\

As for lower bounds, no general bounds are possible for $A(a,t)$ if $t> \frac 1 2$ since then e.g. $A(e_1,t) = 0$. However, for $t=\frac 1 2$, a lower bound independent from $n \in \N$ and $a \in S^{n-1}$ was shown by K\"onig and Rudelson \cite{KR}, namely
$$\frac 1 {17} \le A(a, \frac 1 2) \le 1 \ . $$
For $n=2, 3$, the minimum is of $A( . , \frac 1 2)$ is attained for $a = a^{[n]}$, cf. K\"onig and Koldobsky \cite{KK1}. Whether this holds for $n>3$, too, is an open question. In \cite{KR} similar lower bounds are also given for the volume of $k$-codimensional sections of the cube at distance $t=\frac 1 2$ from the origin. They are of the form $c^k$ where $c$ is an absolute constant. \\
The results of \cite{KR} and part (b) of Proposition \ref{prop4.3} allow a similar lower bound for the perimeter.

\begin{proposition}\label{prop4.4}
Let $n \ge 3$, $K=Q^n$ and $a \in S^{n-1}$. Then
$$\frac n {17} \le P(a,\frac 1 2) \le P(a,0) = 2 ( (n-2) \sqrt 2 + 1 ) \ . $$
\end{proposition}

\begin{proof}
Let $a \in S^{n-1}$, $a_1 > 0$ and put $a_j' := \frac{a_j}{\sqrt{1-a_1^2}}$, $j=1 \etc n$, $\tilde{a}' := (a_j')_{j=2}^n \in S^{n-2} \subset \R^{n-1}$. Then
\begin{align*}
\sqrt{1-a_1^2} & \ \frac 2 {\pi} \int_0^\infty \prod_{j=2}^n \frac{\sin(a_j r)}{a_j r} \ \cos(a_1 r) \ \cos(r) \ dr \\
& = \frac 2 {\pi} \int_0^\infty \prod_{j=2}^n \frac{\sin(a_j' s)}{a_j' s} \ \cos(a_1' s) \ \cos(\frac s {\sqrt{1-a_1^2}}) \ ds \\
& = \frac 1 {\pi} \int_0^\infty \prod_{j=2}^n \frac{\sin(a_j' s)}{a_j' s} \ [ \ \cos(\frac s {\sqrt{1-a_1}}) + \cos(\frac s {\sqrt{1+a_1}}) \ ] \ ds \\
& = \frac 1 2 ( A_{(n-2)} (\tilde{a}', \frac 1 {2 \sqrt{1-a_1}}) + A_{(n-2)} (\tilde{a}', \frac 1 {2 \sqrt{1+a_1}}) ) \ ,
\end{align*}
where $ A_{(n-2)} (\tilde{a}', \frac 1 {2 \sqrt{1 \mp a_1}})$ denotes the $(n-2)$-dimensional volume of the section of the $(n-1)$-cube by a hyperplane orthogonal to $\tilde{a}'$ at distance
$\frac 1 {2 \sqrt{1 \mp a_1}}$. Since $\frac 1 {2 \sqrt{1 + a_1}} \le \frac 1 2$, we have
$A_{(n-2)} (\tilde{a}', \frac 1 {2 \sqrt{1+a_1}}) \ge A_{(n-2)}(\tilde{a}',\frac 1 2) \ge \frac 1 {17}$. Summing up these estimates, we find using Proposition \ref{prop4.3} that $P(a,\frac 1 2) \ge \frac 1 {17}$. Thus $P(a,t)$ is of order $n$ for any $0 \le t \le \frac 1 2$.
\end{proof}


\end{document}